\documentclass[a4paper]{amsart}

\usepackage{amssymb}
\usepackage{amsmath,amsthm}
\usepackage{enumerate,paralist}
\usepackage{amstext}
\usepackage{psfrag}
\usepackage{dsfont}
\usepackage[dvips]{epsfig}
\usepackage[english]{babel}
\usepackage{color}

\theoremstyle{plain}
\newtheorem{theorem}{Theorem}
\newtheorem{corollary}{Corollary}

\theoremstyle{definition}

\newtheorem{example}{Example}
\newtheorem{remark}{Remark}
\newtheorem{assumption}{Assumption}

\numberwithin{theorem}{section}
\numberwithin{corollary}{section}
\numberwithin{lemma}{section}
\numberwithin{definition}{section}
\numberwithin{example}{section}
\numberwithin{remark}{section}
\numberwithin{proposition}{section}
\numberwithin{assumption}{section}

\renewcommand{\leq}{\leqslant}
\renewcommand{\geq}{\geqslant}
\newcommand{\intl}{\int\limits}
\newcommand{\liml}{\lim\limits}
\newcommand{\suml}{\sum\limits}

\newcommand{\cC}{\mathbb C}

\newcommand{\cR}{\mathbb R}

\newcommand{\cRd}{{{\mathbb R}^d}}

\newcommand{\Nat}{\mathbb N}

\newcommand{\E}{\mathbb E}

\newcommand{\pd}{\partial}
\newcommand{\ffi}{\varphi}
\newcommand{\eps}{\varepsilon}
\newcommand{\any}{\forall}

\newcommand{\sbs}{\subset}
\newcommand{\supp}{\mathop{\mathrm{supp}}\nolimits}

\newcommand{\id}{\mathop{\mathrm{Id}}\nolimits}
\newcommand{\Id}{\mathop{\mathrm{Id}}\nolimits}
\newcommand{\Dom}{\mathop{\mathrm{Dom}}\nolimits}
\newcommand{\D}{\mathop{\mathrm{Dom}}\nolimits}

\newcommand{\Hess}{\mathop{\mathrm{Hess}}\nolimits}
\newcommand{\tr}{\mathop{\mathrm{tr}}\nolimits}

\newcommand{\bu}{\bullet}

\newcommand{\Sharp}[1]{\wideparen{#1}}
\usepackage{MnSymbol}
\newcommand{\RE}{\mathop{\mathrm{Re}\,}\nolimits}
\begin{document}

\title[The method of Chernoff approximation]{The method of Chernoff approximation}

\author{ Yana A. Butko }

\address{Saarland University, Postfach 151150, D-66041, Saarbr\"{u}cken, Germany }
\date{\today}

\begin{abstract}
This survey describes the method of approximation of operator semigroups, based on the Chernoff theorem. We outline recent results in this domain as well as clarify relations between constructed approximations, stochastic processes, numerical schemes for PDEs and SDEs, path integrals. We discuss Chernoff approximations for operator semigroups and Schr\"{o}dinger groups. In particular, we consider Feller semigroups in $\cRd$,   (semi)groups
obtained from some original (semi)groups by different procedures: additive perturbations of generators, multiplicative perturbations of generators (which sometimes corresponds to a random time-change of related stochastic processes),  
subordination of  semigroups / processes, imposing boundary / external conditions (e.g., Dirichlet or Robin conditions), averaging of generators, ``rotation'' of semigroups. 
The developed  techniques can be combined to approximate (semi)groups  obtained via several iterative procedures listed above. Moreover, this method can be implemented to obtain approximations for solutions of some time-fractional evolution equations, although these solutions do not posess the semigroup property.

\bigskip

\end{abstract}

\maketitle

\tableofcontents

\section{Introduction}
Let $(X,\|\cdot\|_X)$ be a Banach space. A family $(T_t)_{t\geq0}$ of bounded linear operators on $X$ is called a \emph{strongly continuous semigroup} (denoted as \emph{$C_0$-semigroup}) if $T_0=\id$, $T_t\circ T_s=T_{t+s}$ for all $t,s\geq0$, and $\lim_{t\to0}\|T_t\ffi-\ffi\|_X=0$ for all $\ffi\in X$. The \emph{generator} of the semigroup $(T_t)_{t\geq0}$ is an operator $(L,\Dom(L))$ in $X$ which is given by $L\ffi:=\lim_{t\to0}t^{-1}(T_t\ffi-\ffi)$, $\Dom(L):=\left\{\ffi\in X\,:\, \lim_{t\to0}t^{-1}(T_t\ffi-\ffi)\,\text{ exists in }\, X \right\}$. In the sequel, we denote the semigroup with a given generator $L$ both as  $(T_t)_{t\geq0}$ and as $(e^{tL})_{t\geq0}$. The following fundamental result of the theory of operator semigroups  connects $C_0$-semigroups and evolution equations: \emph{let $(L,\Dom(L))$ be a densely defined linear operator in $X$ with a nonempty resolvent set. The Cauchy problem $\frac{\pd f}{\pd t}=Lf$, $f(0)=f_0$ in $X$ for every $f_0\in\Dom(L)$ has a unique solution $f(t)$ which is continuously differentiable on $[0,+\infty)$   if and only if $(L,\Dom(L))$ is the generator of a $C_0$-semigroup $(T_t)_{t\geq0}$ on $X$. And the solution is given by $f(t):=T_tf_0$.} 

Let now $Q$ be a locally compact metric space.  Let $(\xi_t)_{t\geq0}$ be a temporally homogeneous Markov process with the state space $Q$ and with transition probability $P(t,x,dy)$. The family $(T_t)_{t\geq0}$, given by $T_t\ffi:=\int_Q \ffi(y) P(t,x,dy)$, is a semigroup which, for several important classes of Markov processes, happens to be strongly continuous on some suitable Banach spaces of functions on $Q$. Hence, in this case, we have three equivalent problems:

\begin{enumerate}
\item to construct the $C_0$-semigroup $(T_t)_{t\geq0}$ with a given generator $(L,\Dom(L))$ on a given Banach space $X$;

\item to solve the Cauchy problem $\frac{\pd f}{\pd t}=Lf$, $f(0)=f_0$ in $X$;

\item to determine the transition kernel $P(t,x,dy)$ of an underlying Markov process $(\xi_t)_{t\geq0}$.
\end{enumerate}
The basic example is given by the operator $(L,\Dom(L))$ which is the closure of $(\frac12\Delta,S(\cRd))$ in the Banach space\footnote{We denote the space of continuous functions on $\cRd$ vanishing at infinity by $C_\infty(\cRd)$ and the Schwartz space by  $S(\cRd)$.} $X=C_\infty(\cRd)$ or in $X=L^p(\cRd)$, $p\in[1,\infty)$.  The operator $(L,\Dom(L))$ generates a $C_0$-semigroup  $(T_t)_{t\geq0}$ on $X$; this semigroup is given for each $f_0\in X$ by
\begin{align}\label{heatSG}
T_tf_0(x)=(2\pi t)^{-d/2}\intl_{\cRd}f_0(y)\exp\left\{-\frac{|x-y|^2}{2t} \right\}dy;
\end{align}
the function $f(t,x):=T_tf_0(x)$ solves the corresponding Cauchy problem for the heat equation $\frac{\pd f}{\pd t}=\frac12\Delta f$; and 
\begin{align}\label{heat-P}
P(t,x,dy):=(2\pi t)^{-d/2}\exp\left\{-\frac{|x-y|^2}{2t} \right\}dy
\end{align}
 is the transition probability of a $d$-dimensional Brownian motion. However, it is usually not possible to determine a $C_0$-semigroup in an explicit form, and one has to approximate it. In this note, we demonstrate the method of approximation based on the Chernoff theorem (\cite{MR0231238,MR0417851}). In the sequel, we use the following (simplified) version of the Chernoff theorem, assuming that the existence of the semigroup under consideration is already established.
\begin{theorem}\label{ChernoffTheorem}
Let $(F(t))_{t\geq0}$ be a family of bounded linear operators on a Banach space $X$. Assume that
\begin{enumerate}
\item[(i)]  $F(0)=\id$,

\item[(ii)] $\|F(t)\|\leq  e^{wt}  $ for some   $w\in\cR$  and all $t\geq0$,

\item[(iii)] the limit $L\ffi:=\liml_{t\to0}\frac{F(t)\ffi-\ffi}{t}$ exists for all $\ffi\in D$, where $D$ is a dense subspace in $X$ such that $(L,D)$ is closable and the closure $(L,\Dom(L))$ of $(L,D)$ generates a $C_0$-semigroup  $(T_t)_{t\geq0}$.
\end{enumerate}
 Then the semigroup $(T_t)_{t\geq0}$ is  given by
\begin{equation}\label{0:eq:ChernoffFormula}
T_t\ffi=\liml_{n\to\infty}[F(t/n)]^n\ffi
\end{equation}
for all $\ffi\in X$, and the convergence is locally uniform with respect to  $t\geq0$. 
\end{theorem}
Any family $(F(t))_{t\geq0}$ satisfying the assumptions of the Chernoff theorem~\ref{ChernoffTheorem} with respect to a given $C_0$-semigroup $(T_t)_{t\geq0}$ is called \emph{Chernoff equivalent}, or \emph{Chernoff tangential}  to the semigroup $(T_t)_{t\geq0}$. And the formula~\eqref{0:eq:ChernoffFormula} is  called \emph{Chernoff approximation} of $(T_t)_{t\geq0}$.
 Evidently, in the case of a bounded generator $L$, the family $F(t):=\id+tL$ is Chernoff equivalent to the semigroup $(e^{tL})_{t\geq0}$. And we get a classical formula 
 \begin{align}\label{formula-E}
 e^{tL}=\lim_{n\to\infty}\left[\id+\frac{t}{n}L \right]^n.
 \end{align}
  Moreover, for an arbitrary generator $L$, one considers $F(t):=\left(\Id-tL\right)^{-1}\equiv \frac{1}{t}R_L(1/t)$ and obtains the \emph{Post--Widder inversion formula}:
$$
T_t \ffi=\liml_{n\to\infty}\left(\Id-\frac{t}{n}L\right)^{-n}\ffi\equiv \liml_{n\to\infty}\left[\frac{n}{t}R_L(n/t)\right]^n\ffi,\quad \forall \,\ffi\in X.
$$
A well-developed functional calculus approach to Chernoff approximation of  $C_0$-semigroups by families $(F(t))_{t\geq0}$, which are given by (bounded completely monotone) functions of the generators (as, e.g.,  in the case of the Post--Widder inversion formula above), can be found in~\cite{Tomilov}. We use another approach. We are looking for arbitrary families  $(F(t))_{t\geq0}$ which are Chernoff equivalent to a given $C_0$-semigroup (i.e., the only connection of $F(t)$ to  the generator $L$ is given  via the assertion (iii) of the Chernoff theorem). But we are especially interested in families $(F(t))_{t\geq0}$ which are given explicitly (e.g., as integral operators with explicit kernels or pseudo-differential operators with explicit symbols). This is useful both for practical calculations and for further interpretations of Chernoff approximations as path integrals (see, e.g., \cite{MR2863557,MR2423533,MR3455669} and references therein). Moreover, we consider different operations on generators (what sometimes corresponds to operations  on Markov processes) and  find out, how to construct Chernoff approximations for $C_0$-semigroups with modified generators on the base of Chernoff approximations for original ones.
\begin{center}
{\scalebox{0.25}{\includegraphics{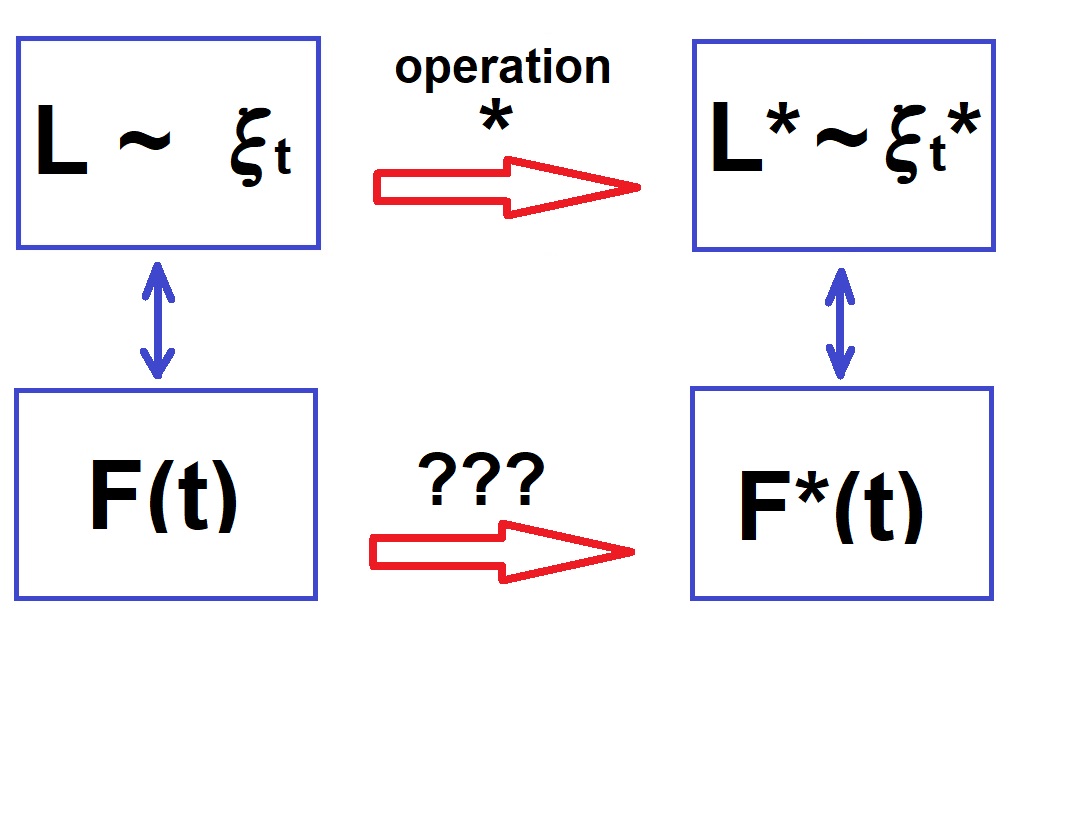}}}
\end{center}

\vspace{-0.6cm}
This approach allows to create a kind of a LEGO-constructor: we start with a $C_0$-semigroup which is already known\footnote{E.g., the semigroup generated by a Brownian motion on a star graph with Wentzell boundary conditions at the vertex~\cite{MR2927703}; see also~\cite{MR1912205,MR3231629,MR2042661} for further examples.} or Chernoff approximated\footnote{E.g., the  semigroup generated by a Brownian motion on a compact Riemannian manifold (see~\cite{MR2276523} and references therein) and (Feller) semigroups generated by Feller processes in $\cRd$ (see~\cite{MR2999096} and Sec.~\ref{subsec_Feller}).}; then, applying different operations on its generator, we consider more and more complicated $C_0$-semigroups and  construct their Chernoff approximations.
\begin{center}
{\scalebox{0.25}{\includegraphics{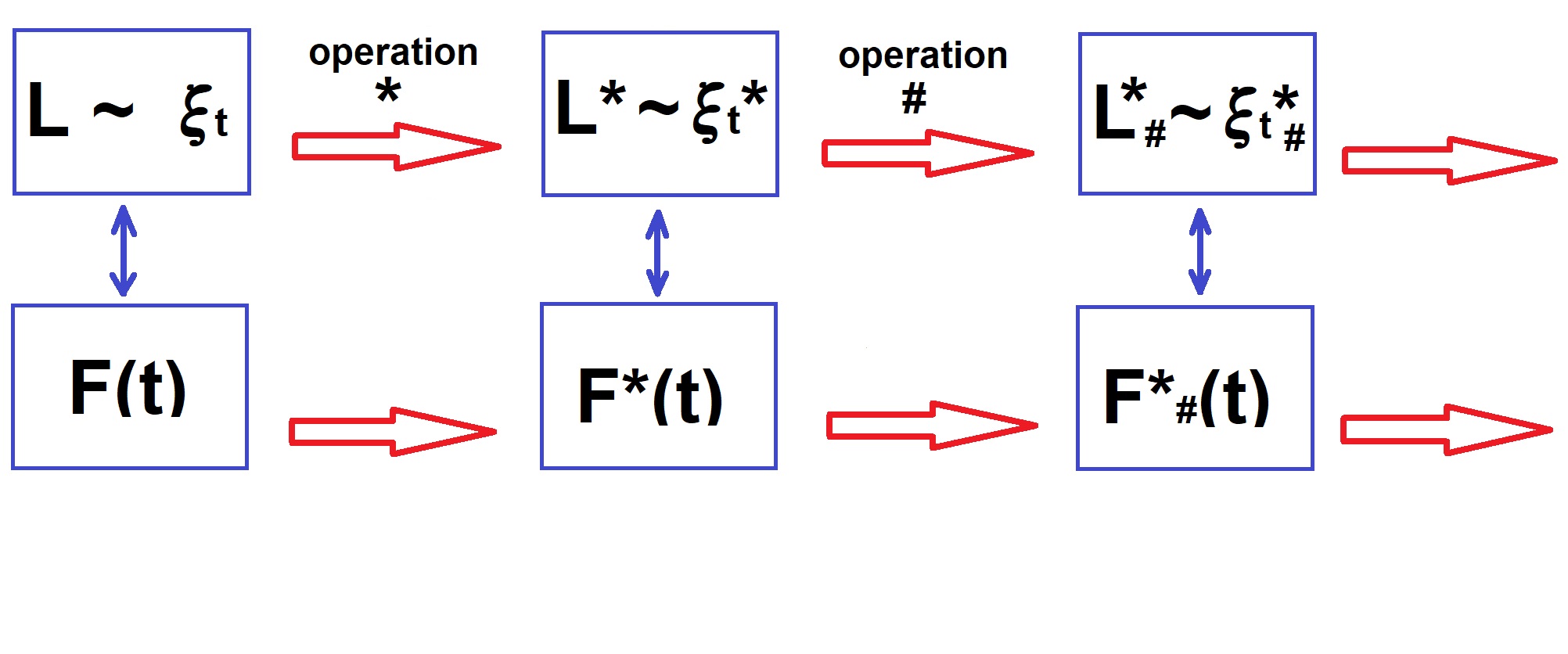}}}
\end{center} 

\vspace{-0.6cm}
\noindent Chernoff approximations are available for the following operations:
%\begin{itemize}

\noindent $\bu$ Operator splitting; additive perturbations of a generator (Sec.~\ref{subsec_additive}, \cite{Habil,technomag,MR2999096});

\noindent  $\bu$ Multiplicative perturbations of a generator / random time change of a process via an additive functional (Sec.~\ref{subsec_mult}, \cite{Habil,technomag,MR2999096});

\noindent  $\bu$ killing of a process upon leaving a given domain / imposing Dirichlet boundary (or external) conditions (Sec.~\ref{subsec_BC}, \cite{Habil,MR3903609,MR2729591});

\noindent $\bu$  imposing  Robin  boundary conditions (Sec.~\ref{subsec_BC}, \cite{MR2541275});

\noindent $\bu$ subordination of a semigroup / process (Sec.~\ref{subsec_subordinate}, \cite{Habil,MR3804267});

\noindent $\bu$ ``rotation'' of a semigroup (see Sec.~\ref{subsec_Schr}, \cite{MR3490776,MR3871530});

\noindent $\bu$ averaging of  semigroups (see Sec.~\ref{subsec_Schr}, \cite{MR3819810,MR3479996,MR3819810,MR3844136});

%\end{itemize}
\smallskip

\noindent Moreover, Chernoff  approximations have been obtained for some stochastic Schr\"{o}-dinger type equations  in~\cite{MR3527030,MR2314121,MR2157588,MR2190085}; for evolution equations with the Vladimirov operator (this operator  is a $p$-adic analogue of the Laplace operator)  in~\cite{MR2963683,MR2865952,MR2681357,MR2599557,MR2462100}; for evolution equations containing L\'{e}vy Laplacians in~\cite{MR2359376,MR2348706}; for some nonlinear equations in~\cite{MR3791362}.

%\subsection{Interpretations of Chernoff approximation}
%or Remark
%!!!!!!!!!!!!!!!!!!!!!!!!!1

 Chernoff approximation can be interpreted  as a numerical scheme for solving evolution equations. Namely, for the Cauchy problem  $\frac{\pd f}{\pd t}=Lf$, $f(0)=f_0$, we have:
 $$
  u_0:=f_0,\quad\quad u_{k}:=F(t/n)u_{k-1},\quad k=1,\ldots,n,\quad \quad \quad f(t)\approx u_n.
$$
In some particular cases, Chernoff approximations are an abstract analoge of the \emph{operator splitting method}  known in the numerics of PDEs (see Remark~\ref{1:rem:F-theta}). And the Chernoff theorem itself can be understood as a version of the \emph{``Meta-theorem of numerics'': consistency and stability imply convergence.} Indeed, conditions (i) and (iii) of Theorem~\ref{ChernoffTheorem} are consistency conditions, whereas condition (ii) is a stability condition. Moreover, in some cases, the families $(F(t))_{t\geq 0}$ give rise to Markov chain approximations  for  $(\xi_t)_{t\geq0}$ and provide Euler--Maruyama schemes for the corresponding SDEs (see Example~\ref{HFF-Gaussian}).

 If all operators $F(t)$ are integral operators with elementary kernels or pseudo-differential operators with elementary symbols, the identity  \eqref{0:eq:ChernoffFormula} leads to representation of a given semigroup by $n$-folds iterated integrals of elementary functions when $n$ tends to infinity.  This gives rise to \emph{Feynman formulae}.
%\begin{definition}
A \emph{Feynman formula} is a representation of a solution of an initial (or initial-boundary) value problem for an evolution equation (or, equivalently, a representation of the semigroup solving the problem) by a limit of $n$-fold iterated integrals of some functions as $n\to\infty$.
%\end{definition}
%\begin{remark}
One should not confuse the notions of Chernoff approximation and Feynman formula.
On the one hand, not all Chernoff approximations can be directly interpreted as Feynman formulae since, generally, the operators $(F(t))_{t\geq0}$ do not have to be neither integral operators, nor pseudo-differential operators.  On the other hand,   representations of solutions of evolution equations in the form of Feynman formulae can be obtained by different methods, not necessarily via the Chernoff Theorem. And such Feynman formulae may have no relations to any  Chernoff approximation, or their relations may be quite indirect.
Richard Feynman was the first who considered representations of solutions of evolution equations by limits of iterated integrals  (\cite{MR0026940}, \cite{MR0044379}).  He has, namely,  introduced  a construction of a path integral (known nowadays as \emph{Feynman path integral})  for solving the Schr\"{o}dinger equation. And this path integral was defined exactly as a  limit of  iterated finite dimensional integrals. Feynman path integrals can be also understood as integrals with respect to  Feynman type pseudomeasures. Analogously, one can sometimes obtain representations of  a solution of an initial (or initial-boundary) value problem for an evolution equation (or, equivalently, a representation of an operator semigroup resolving the problem) by functional (or, path) integrals with respect to probability measures. Such representations are usually   called  \emph{Feyn\-man--Kac formulae}. It is a usual situation that  limits in  Feynman formulae coincide with (or in some cases define) certain path integrals with respect to probability measures or Feynman type pseudomeasures on a set  of paths of a physical system.  Hence the iterated integrals in   Feynman formulae for some problem give approximations to  path integrals  representing the solution of the same problem.  Therefore,  representations of evolution semigroups by Feynman formulae, on the one hand, allow to establish new path-integral-representations and, on the other hand, provide an additional tool to calculate path integrals numerically. Note  that different  Feynman formulae  for the same semigroup allow to establish relations between different path integrals (see, e.g., \cite{Habil}).

%\begin{remark}
The result of Chernoff  has diverse generalizations. Versions, using  arbitrary partitions  of the time interval $[0,t]$ instead of the equipartition $(t_k)_{k=0}^n$ with $t_k-t_{k-1}=t/n$, are presented, e.g., in \cite{MR710486}, \cite{MR2013201}. The analogue of the Chernoff theorem for  multivalued generators can be found, e.g., in  \cite{MR838085}. Analogues  of Chernoff's result for  semigroups, which are continuous in a weaker sense, are obtained, e.g., in \cite{MR2026919}, \cite{DisserKuehnemund}. For analogues of the Chernoff theorem in the case of nonlinear semigroups, see, e.g.,  \cite{MR0390843}, \cite{MR0448185}, \cite{MR0293452}. The Chernoff Theorem for two-parameter families of operators  can be found in \cite{MR2157588}, \cite{MR2965550}.

%%\cite{MR710486} Theo.~3.5.3; \cite{MR1721989} Theo.~III.5.2; \cite{JacobI} Theo. 4.4.13; 
%\end{remark}

\section{Chernoff approximations for operator semigroups and further applications}
\subsection{Chernoff approximations for the procedure of operator splitting}\label{subsec_additive}
\begin{theorem}%[Chernoff approximation of semigroups generated by additive perturbations of generators]
\label{1:cor:additivePert}
Let $(T_t)_{t\geq0}$ be a strongly continuous semigroup on a Banach space $X$ with generator $(L,\Dom(L))$. Let $D$ be a core for $L$. Let  $L=L_1+\ldots+L_m$ hold on $D$ for some linear operators $L_k$, $k=1,\ldots,m$, in  $X$. %with $D\subset\Dom(L_k)$ for all $k=1,\ldots ,m$.
  Let  $(F_k(t))_{t\geq0}$, $k=1,\ldots ,m,$ be families of bounded linear operators on $X$  such that  for all   $k\in\{1,\ldots ,m \}$ holds: $F_k(0)=\id$,  $\|F_k(t)\|\leq e^{a_kt}$ for some $a_k>0$  and all $t\geq0$, %as well as 
    $\lim_{t\to0}\big\|
\frac{F_k(t)\ffi-\ffi}{t}-L_k\ffi\big\|_X=0$ for all  $\ffi\in D$. Then the family  $(F(t))_{t\geq0}$, with $F(t):= F_1(t)\circ\cdots \circ F_m(t)$,
is Chernoff equivalent to the semigroup $(T(t))_{t\geq0}$. And hence the Chernoff approximation
\begin{align}\label{ChAppOperatorSplitting}
T_t\ffi=\liml_{n\to\infty} \big[F(t/n)   \big]^n\ffi\equiv \liml_{n\to\infty} \big[F_1(t/n)\circ\cdots \circ F_m(t/n)   \big]^n\ffi
\end{align}
holds for each $\ffi\in X$ locally uniformly with respect to $t\geq 0$. 
\end{theorem}
Note that we do not require from summands $L_k$ to be generators of $C_0$-semigroups. For example, $L_1$ can be a leading term (which generates a $C_0$-semigroup) and $L_2,\ldots,L_m$ can be $L_1$-bounded additive perturbations such  that  $L:=L_1+L_2+\cdots+L_m$ again generates a strongly continuous semigroup. Or even $L$ can be a sum of operators $L_k$, none of which generates a strongly continuous semigroup itself. %The proof of this theorem is  the same as the proof of Theorem~5.1 of~\cite{MR2999096}.
%Obviously, the family $(F(t))_{t\geq0}$ is strongly continuous if and only if all the families $(F_k(t))_{t\geq0}$, $k=1,\ldots ,m,$ are strongly continuous.
\begin{proof}%%%\noindent {\bfseries Proof:}
%Since the families $(F_k(t))_{t\geq0}$, $k=1,\ldots ,m,$  are Chernoff equivalent to the semigroups $(T_k(t))_{t\geq0}$ respectively,  we have $F_k(0)=\id$ and  $\|F_k(t)\|\leq e^{a_kt}$ for some $a_k>0$  for each $k\in\{1,\ldots ,m \}$.
 Obviously, the family $(F(t))_{t\geq0}$ satisfies  the conditions   $F(0)=\id$ and
$\|F(t)\|\leq \|F_1(t)\|\cdot\ldots \cdot\|F_m(t)\|\leq e^{(a_1+\cdots +a_m)t}.$ 
Further, for each $\ffi\in D$, we have
\begin{align*}
&\lim_{t\to0}\bigg\|
\frac{F(t)\ffi-\ffi}{t}-L\ffi\bigg\|_X
=\lim_{t\to0}\bigg\|
\frac{F_1(t)\circ\cdots \circ F_m(t)\ffi-\ffi}{t}-L_1\ffi-\cdots -L_m\ffi\bigg\|_X\\
&=
\lim_{t\to0}\bigg\|
F_1(t)\circ\cdots \circ F_{m-1}(t)\left(\frac{F_m(t)\ffi-\ffi}{t}- L_m\ffi\right)\\
&
+\left(F_1(t)\circ\cdots \circ F_{m-1}(t)-\Id\right)L_m\ffi+\frac{F_1(t)\circ\cdots \circ F_{m-1}(t)\ffi-\ffi}{t}- L_1\ffi-\cdots -L_{m-1}\ffi   \bigg\|_X\\
&
 \leq
\lim_{t\to0}\bigg\|
\frac{F_1(t)\circ\cdots \circ F_{m-1}(t)\ffi-\ffi}{t}-L_1\ffi-\cdots -L_{m-1}\ffi\bigg\|_X\\
& 
\leq \cdots \leq
\lim_{t\to0}\bigg\|
\frac{F_1(t)\ffi-\ffi}{t}-L_1\ffi\bigg\|_X
=0.
%\qedhere
\end{align*}
Therefore, all requirements of the Chernoff theorem~\ref{ChernoffTheorem} are fulfilled and hence $(F(t))_{t\geq0}$ is Chernoff equivalent to $(T(t))_{t\geq0}$.
\end{proof}

\begin{remark}\label{1:rem:F-theta}
Let all the assumptions of Theorem~\ref{1:cor:additivePert} be fulfilled. Consider for simplicity the case $m=2$.  Let $\theta,\tau\in[0,1]$. Similarly to the proof of Theorem~\ref{1:cor:additivePert}, one shows that the following families $(H^\theta(t))_{t\geq0}$ and $(G^{\tau}(t))_{t\geq0}$ are Chernoff equivalent to the semigroup $(T_t)_{t\geq0}$ generated by $L=L_1+L_2$:
\begin{align*}
& H^\theta(t):=F_1(\theta t)\circ F_2(t)\circ F_1((1-\theta)t),\\
&
G^\tau(t):=\tau F_1(t)\circ F_2(t)+ (1-\tau)F_2(t)\circ F_1(t). 
\end{align*}
Note that we have $H^0(t)=F_2(t)\circ F_1(t)$, and  $H^1(t)=F_1(t)\circ F_2(t)$.  Hence the parameter $\theta$ corresponds to different orderings of non-commuting terms $F_1(t)$ and $F_2(t)$.  Further,  $G^{1/2}(t)=\frac12\left(H^1(t)+H^0(t)\right)$. In the case when both $L_1$ and $L_2$ generate $C_0$-semigroups and $F_k(t):=e^{tL_k}$, Chernoff approximation~\eqref{ChAppOperatorSplitting}  with families $(H^\theta(t))_{t\geq0}$, $\theta=1$ or $\theta=0$, reduces to the classical Daletsky--Lie--Trotter formula.  Moreover, Chernoff approximation~\eqref{ChAppOperatorSplitting}  can be understood as  an abstract analogue of the \emph{operator splitting} known in numerical methods of solving PDEs (see \cite{MR3617561} and references therein).   If  $\theta=0$ and $\theta=1$, the families $(H^\theta(t))_{t\geq0}$ correspond to first order splitting schemes. Whereas the family $(H^{1/2}(t))_{t\geq0}$  corresponds to the symmetric Strang splitting and, together with $(G^{1/2}(t))_{t\geq0}$, represents second order splitting schemes.
\end{remark}

%%%%%%%%%%%%%%%%%%%%%%%%%%%%%%%%%%%%%%%%%%%%%%%%%%%%%%%%%%%%%%%%%%%%%%%%%%%%%%%%%%%%%%%%%%%%%%%%%%%%%%%%%%%%%

\subsection{Chernoff approximations for Feller semigroups}\label{subsec_Feller}
We consider the Banach space $X=C_\infty(\cRd)$ of continuous functions on $\cRd$, vanishing at infinity.
A semigroup of bounded linear operators $(T_t)_{t\geq0}$ on the Banach space $X$ is called \emph{Feller semigroup} if it is a strongly continuous semigroup, it is \emph{positivity preserving} (i.e. $T_t\ffi\geq0$ for all $\ffi\in X$ with $\ffi\geq0$) and it is \emph{sub-Markovian} (i.e. $T_t\ffi\leq1$ for all $\ffi\in X$ with $\ffi\leq1$). A Markov process, whose semigroup is Feller, is called  \emph{Feller process}. Let $(L,\Dom(L))$ be the generator of a Feller semigroup $(T_t)_{t\geq0}$.  Assume that $C^\infty_c(\cRd)\subset\Dom(L)$ (this assumption  is quite standard and holds in many cases, see, e.g., \cite{MR3156646}). Then  we have also\footnote{$C^m_\infty(\cRd):=\{\ffi\in C^m(\cRd)\,:\,\pd^\alpha\ffi\in C_\infty(\cRd),\,|\alpha|\leq2 \}$. }  $C^2_\infty(\cRd)\subset\Dom(L)$. And  $L\ffi(x)$ is given for each $\ffi\in C^2_\infty(\cRd)$ and each $x\in\cRd$  by the following formula:
\begin{equation}\label{2:eq:g35}
\begin{split}
    L\ffi(x)
    &=
    -C(x)\ffi(x) - B(x)\cdot\nabla \ffi(x)
    + \tr(A(x)\Hess\ffi(x))\\
    &\phantom{==}+ \int_{y\neq 0} \left( \ffi(x+y) - \ffi(x)
    - \frac{y\cdot\nabla \ffi(x)}{1+|y|^2}\right)\,N(x,dy),
\end{split}
\end{equation}
where  $\Hess\ffi$ is the Hessian matrix of  second order partial derivatives of $\ffi$; as well as  $C(x)\geq0$,  $B(x)\in\cRd$, $A(x)\in\mathbb{R}^{d\times d}$ is a symmetric positive semidefinite matrix    and $N(x,\cdot)$ is a Radon measure  on $\cRd\setminus \{0\}$ with  $\int_{y\neq 0}
|y|^2(1+|y|^2)^{-1}\,N(x,dy)<\infty$  for each $x\in\cRd$.
Therefore,  $L$ is an integro-differential operator on $C^2_\infty(\cRd)$ which is non-local if $N\neq0$.
This class of generators  $L$ includes, in particular,  fractional Laplacians $L=-(-\Delta)^{\alpha/2}$ and relativistic Hamiltonians $\sqrt[\alpha]{(-\Delta)^{\alpha/2}+m(x)}$, $\alpha\in(0,2)$, $m>0$.
Note that  the restriction of $L$ onto $C^\infty_c(\cRd)$ is given by a pseudo-differential operator (PDO)
\begin{align}\label{PDO_L}
L\ffi(x):=(2\pi)^{-d}\intl_{\cRd}\intl_{\cRd}e^{ip\cdot(x-q)}H(x,p)\ffi(q)\,dq\,dp, \quad x\in\cR^d,
\end{align}
with the symbol $-H$ such that
\begin{align}\label{2:eq:g116}
H(x,p)= C(x) +i B(x)\cdot p + p\cdot A(x)p
      + \!\!\!\intl_{y\neq 0}\!\!\!\!
        \left(1-e^{iy\cdot p} + \frac{iy\cdot p}{1+|y|^2} \right)\!\!N(x,dy).
\end{align}
If the symbol $H$ does not depend on $x$, i.e. $H=H(p)$, then the semigroup  $(T_t)_{t\geq0}$ generated by  $(L,\Dom(L))$ is given by (extensions of) PDOs with symbols $e^{-tH(p)}$:
\begin{align*}
T_t\ffi(x)=(2\pi)^{-d}\intl_{\cRd}\intl_{\cRd}e^{ip\cdot(x-q)}e^{-tH(p)}\ffi(q)\,dq\,dp, \quad x\in\cR^d,\quad\ffi\in C^\infty_c(\cRd).
\end{align*}
If the symbol $H$ depends on both variables $x$ and $p$ then  $(T_t)_{t\geq0}$ are again PDOs. However their symbols do not coincide with $e^{-tH(x,p)}$ and are not known explicitly. The family $(F(t))_{t\geq0}$ of PDOs with symbols $e^{-tH(x,p)}$ is not a semigroup any more. However, this family is Chernoff equivalent to $(T_t)_{t\geq0}$. Namely, the following theorem holds (see \cite{MR2999096,MR2759261}):

\begin{theorem}\label{2:thm:HFF}
Let   $H:\cRd\times\cRd\to\cC$ be  measurable, locally bounded in both variables $(x,p)$,  satisfy for each fixed $x\in\cRd$ the  representation \eqref{2:eq:g116} and the following assumptions:
\begin{align*}
&(i)\quad\displaystyle\sup_{q\in\cRd} |H(q,p)| \leq
                 \kappa(1+|p|^2)\quad \text{for all}\quad p\in\cRd \quad\text{and some}\quad \kappa>0,\\
&
(ii)\quad  \displaystyle p\mapsto H(q,p)\quad \text{is uniformly (w.r.t.} \quad q\in\cRd \text{)
                 continuous at}\quad p=0,\\
&
(iii)\quad \displaystyle q\mapsto H(q,p)\quad \text{is continuous for
all}\quad p\in\cRd.
\end{align*}
 Assume that  the function $H(x,p)$ is such that  the PDO with symbol $-H$ defined on  $C_c^\infty(\cRd)$ is closable and the closure (denoted by $(L,\Dom(L))$)   generates a strongly continuous semigroup $(T_t)_{t\geq0}$ on $X=C_\infty(\cRd)$. Consider now for each $t\geq0$  the PDO $F(t)$ with the symbol  $e^{-tH(x,p)}$, i.e. for  $\ffi\in C_c^\infty(\cRd)$
\begin{align}\label{F(t)-PDO}
F(t)\ffi(x)&=(2\pi)^{-d}\int_{\cRd}\int_{\cRd}e^{ip\cdot (x-q)}e^{-tH(x,p)}\ffi(q)dqdp.
\end{align}
Then the  family $(F(t))_{t\geq0}$ extends to a strongly continuous family on $X$  and is Chernoff equivalent to the semigroup $(T_t)_{t\geq0}$.
\end{theorem}

Note that the extensions of $F(t)$ are given (via integration with respect to $p$ in~\eqref{F(t)-PDO}) by integral operators:
\begin{align}\label{eq:F(t)-NonLocCP}
F(t)\ffi(x)=
\int_{\cRd}\ffi(y)\nu^x_t(dy),
\end{align}
where, for each $x\in\cRd$ and each $t\geq0$, the sub-probability measures $\nu^x_t$  are given via their Fourier transform  $\mathcal{F}\left[\nu_t^x\right](p)=(2\pi)^{-d/2}e^{-tH(x,-p)-ip\cdot x}$.

\begin{example}\label{HFF-Gaussian}
Let    in formula~\eqref{2:eq:g35} additionally  $N(x,dy)\equiv0$, the coefficients $A$, $B$, $C$  be bounded and continuous, and  
\begin{align}\label{eq:CP:uniformEllipt}
&\text{there exist } a_0, A_0\in\cR \text{ with } 0<a_0\leq A_0<\infty \text{ such that }\nonumber\\
& a_0|z|^2\leq z\cdot A(x)z\leq A_0|z|^2\quad\text{ for all }\,\,x,z\in\cRd.
\end{align}
Then $L$ is a second order uniformly elliptic operator and  the family $(F(t))_{t\geq0}$  in~\eqref{eq:F(t)-NonLocCP}
has the following view:   $F(0):=\Id$ and for all $t>0$ and all $\ffi\in X$
  \begin{align}\label{F(t)-Gaus}
F(t)\ffi(x):=\frac{e^{-tC(x)}}{\sqrt{(4\pi t)^{d}\det A(x)}}\intl_{\cRd}
e^{-\frac{A^{-1}(x)(x-tB(x)-y)\cdot(x-tB(x)-y)}{4t}}\ffi(y)dy.
\end{align}
Moreover, it has been shown in~\cite{MR2729591} that $F'(0)=L$ on a bigger core $C^{2,\alpha}_c(\cRd)$ what is important for further applications (e.g., in Sec.~\ref{subsec_BC}).

Let now $C\equiv 0$. The evolution equation 
$$
\frac{\pd f}{\pd t}(t,x)=- B(x)\cdot\nabla f(t,x)
    + \tr(A(x)\Hess f(t,x))
$$
is the backward Kolmogorov equation for a $d$-dimesional It\^{o} diffusion process $(\xi_t)_{t\geq0}$ satisfying the SDE
\begin{align}\label{SDE}
d\xi_t=-B(\xi_t)dt+\sqrt{2A(\xi_t)}dW_t,
\end{align}
with a $d$-dimensional Wiener process $(W_t)_{t\geq0}$. Consider the Euler--Maruyama scheme for the SDE~\eqref{SDE} on $[0,t]$ with time step $t/n$:
\begin{align}\label{EulerScheme}
X_0:=\xi_0,\qquad X_{k+1}:=X_k-B(X_k)\frac{t}{n}+\sqrt{\frac{2t}{n}A(X_k)}Z_k,\quad k=0,\ldots,n-1,
\end{align}
where $(Z_k)_{k=0,\ldots,n-1}$ are i.i.d. $d$-dimensional $N(0,\id)$ Gaussian random variables such that $X_k$ and $Z_k$ are independent for all $k=0,\ldots,n-1$. Then, for all $k=0,\ldots,n-1$ holds:
\begin{align*}
\mathbb{E}[f_0(X_{k+1})\,|\,X_k]=\left.\mathbb{E}\left[f_0\left(x-B(x)\frac{t}{n}+\sqrt{\frac{2t}{n}A(x)}Z_k\right)\right]\right|_{x:=X_k}=F(t/n)f_0(X_k).
\end{align*}
By the tower property of conditional expectation, one has
\begin{align*}
&\mathbb{E}[f_0(X_{n})\,|\,X_0=x]=\mathbb{E}[\mathbb{E}[f_0(X_{n})\,|\,X_{n-1}]\,|\,X_0=x]=\ldots=\\
&
\quad=\mathbb{E}[\ldots\mathbb{E}[\mathbb{E}[f_0(X_{n})\,|\,X_{n-1}]\,|\,X_{n-2}]\ldots\,|\, X_0=x]=F^n(t/n)f_0(x).
\end{align*}
Hence, by Theorem~\ref{2:thm:HFF}, it holds for all $x\in\cRd$
\begin{align*}
\mathbb{E}[f_0(\xi_t)\,|\,\xi_0=x]=T_tf_0(x)=\liml_{n\to\infty}F^n(t/n)f_0(x)=\liml_{n\to\infty}\mathbb{E}[f_0(X_{n})\,|\,X_0=x].
\end{align*}
And, therefore, the Euler--Maruyama scheme~\eqref{EulerScheme} converges weakly\footnote{The weak convergence of this  Euler--Maruyama scheme is, of course, a classical result, cf.~\cite{MR1214374}.}. The same holds in the general case of Feller processes satisfying assumptions of Theorem~\ref{2:thm:HFF} (see~\cite{MR2847335}). And the corresponding Markov chain approximation  $(X_k)_{k=0,\ldots,n-1}$ of $\xi_t$ consists of increments of L\'{e}vy processes, obtained form the original Feller process by ``freezing the coefficients'' in the generator in a suitable way (see~\cite{MR2502474}).

Let us investigate the family $(F(t))_{t\geq0}$ in~\eqref{F(t)-Gaus} more carefully.
 We have actually
\begin{align}\label{1:eq:F^ABC+}
F(t)\ffi(x)
=e^{-tC(x)}\int_{\cRd}e^{\frac{A^{-1}(x)B(x)\cdot(x-y)}{2}}e^{-t\frac{|A^{-1/2}(x)B(x)|^2}{4}}\ffi(y)p_A(t,x,y)dy,
\end{align}
where %the function $p_A$ is given by
%\begin{equation}\label{1:eq:3:p_A}
$p_A(t,x,y) :=\left((4\pi t)^{d}\det A(x)\right)^{-1/2} \exp\bigg(-\frac{A^{-1}(x)(x-y)\cdot(x-y)}{4t}\bigg).$
%\end{equation}.
Therefore, Theorem~\ref{2:thm:HFF} yields  the following  Feynman formula  for all $t>0$,  $\ffi\in X$ and $x_0\in\cRd$:
\begin{align}\label{1:eq:LFF with F^ABC+}
&T_t\ffi(x_0)=\liml_{n\to\infty}\intl_{\cR^{dn}}e^{-\frac{t}{n}\suml_{k=1}^n C(x_{k-1})}
e^{\frac12\suml_{k=1}^n A^{-1}(x_{k-1})B(x_{k-1})\cdot (x_{k-1}-x_k)} \times\\
&
\times e^{-\frac{t}{4n}\suml_{k=1}^n |A^{-1/2}(x_{k-1})B(x_{k-1})|^2}\ffi(x_n)p_A(t/n,x_0,x_1)\ldots p_A(t/n,x_{n-1},x_n)\,dx_1\cdots dx_n.\nonumber
\end{align}
And the convergence is uniform with respect to $x_0\in\cRd$ and $t\in(0,t^*]$ for all $t^*>0$.
 The limit in the right hand side of  formula~\eqref{1:eq:LFF with F^ABC+} coincides with  the following path integral 
  (compare with the  formula (34) in~\cite{MR1614598} and formula (3) in~\cite{MR2052140}):
\begin{align*}
T_t\ffi(x_0)=&\mathbb{E}^{x_0} \bigg[\exp\bigg(-\intl_0^t C(X_s)ds\bigg)\exp\bigg(-\frac12\intl_0^t A^{-1}(X_s)B(X_s )\cdot dX_s)\bigg) \times\\
&
\phantom{qwqwwwwqw}\times\exp\bigg(-\frac{1}{4}\intl_0^t A^{-1}(X_s)B(X_s)\cdot B(X_s)ds \bigg)\ffi(X_t ) \bigg].\nonumber
\end{align*}
Here the stochastic integral   $\int_0^t A^{-1}(X_s)B(X_s)\cdot dX_s$ is an It\^{o} integral. And $\mathbb{E}^{x_0}$ is the expectation of a (starting at $x_0$) diffusion process  $(X_t)_{t\ge0}$ with the variable diffusion matrix  $A$ and without any  drift, i.e   $(X_t)_{t\ge0}$ solves the stochastic differential equation
$$
dX_t=\sqrt{2A(X_t)}dW_t.
$$
\end{example}

\begin{remark}\label{remarkOnCores}
Let  now  $N(x,dy):=N(dy)$ in formula~\eqref{2:eq:g35}, i.e. $N$ does not depend on $x$. Let  the coefficients $A$, $B$, $C$  be bounded and continuous, and  the property~\eqref{eq:CP:uniformEllipt} hold. Then $L=L_1+L_2$, where $L_1$ is the local part of $L$, given in the first line of ~\eqref{2:eq:g35} and $L_2$ is the non-local part of $L$, given in the second line of~\eqref{2:eq:g35}. And, respectively, $H(x,p)=H_1(x,p)+H_2(p)$ in~\eqref{2:eq:g116}, where $H_1(x,p)$ is a quadratic polynomial with respect to $p$ with variable coefficients and $H_2$ does not depend on $x$. Then the closure of $(L_2, C^\infty_c(\cRd))$ in $X$ generates a $C_0-$semigroup $(e^{tL_2})_{t\geq0}$  and operators  $e^{tL_2}$ are PDOs with symbols $e^{-tH_2}$ on $C^\infty_c(\cRd)$. Let the family of probability measures $(\eta_t)_{t\geq0}$ be such that $\mathcal{F}[\eta_t]=(2\pi)^{-d/2}e^{-tH_2}$. Then we have $e^{tL_2}\ffi=\ffi*\eta_t$ on $X$. Assume that $H_2\in C^\infty(\cRd)$. Then the family $(F(t))_{t\geq0}$ in~\eqref{F(t)-PDO} can be represented (for $\ffi\in C^\infty_c(\cRd)$) in the following way (cf. with formula~\eqref{eq:F(t)-NonLocCP}):
\begin{align*}
F(t)\ffi(x)&=\left[\mathcal{F}^{-1}\circ e^{-tH(x,\cdot)}\circ\mathcal{F}\ffi\right](x)=\left[\mathcal{F}^{-1}\circ e^{-tH_1(x,\cdot)}\circ\mathcal{F}\circ\mathcal{F}^{-1}\circ e^{-tH_2}\circ\mathcal{F}\ffi\right](x)\\
&
=\left(\ffi*\eta_t*\rho^x_t\right)(x),
\end{align*}
where $\rho_t^x(z):=e^{-tC(x)}\left((4\pi t)^{d}\det A(x)\right)^{-1/2}
\exp\left\{-\frac{A^{-1}(x)(z-tB(x))\cdot(z-tB(x))}{4t}\right\},$ i.e. the family $(F_1(t))_{t\geq0}$, $F_1(t)\ffi(x):=(\ffi*\rho^x_t)(x)$, is actually given by formula~\eqref{F(t)-Gaus}. The representation 
\begin{align}\label{F(t)via2conv}
F(t)\ffi(x)=\left(\ffi*\eta_t*\rho^x_t\right)(x)
\end{align}
holds even for all $\ffi\in X$, $x\in\cRd$ and without the assumption that $H_2\in C^\infty(\cRd)$.  Denoting $e^{tL_2}$ as $F_2(t)$, we obtain  that $F(t)=F_1(t)\circ F_2(t)$. Due to Theorem~\ref{2:thm:HFF}, $F'(0)=L$ on a core $D:=C^\infty_c(\cRd)$. Using Theorem~\ref{1:cor:additivePert} and Example~\ref{HFF-Gaussian}, one shows that $F'(0)=L$ even on $D=C^{2,\alpha}_c(\cRd)$ as soon as $C^{2,\alpha}_c(\cRd)\subset\Dom(L_2)$ (without the assumption  $H_2\in C^\infty(\cRd)$). The bigger core $D$ is more suitable for further applications  of the family $(F(t))_{t\geq0}$ in the form of~\eqref{F(t)via2conv} in Sec.~\ref{subsec_BC}. 
\end{remark}

\begin{example} Consider the symbol $H(x,p):=a(x)|p|$, where $a\in C^\infty(\cRd)$ is a strictly positive bounded function.    The closure of the PDO $(L, C^\infty_c(\cRd))$ with symbol $-H$  acts as  $L\ffi(x):=a(x)\left(-(-\Delta)^{1/2} \right)\ffi(x)$, 
generates a Feller semigroup $(T_t)_{t\geq0}$ and, by Theorem~\ref{2:thm:HFF}, the following family $(F(t))_{t\geq0}$ is Chernoff equuivalent to $(T_t)_{t\geq0}$:

\begin{align*}
 F(t)\ffi(x):&=(2\pi)^{-d}\intl_{\cR^{d}}\intl_{\cRd}e^{ip\cdot(x-q)}e^{-{ ta(x)}|p|}\ffi(q)dqdp\\
 &
 =\Gamma\left(\frac{d+1}{2} \right)\intl_{\cRd}\ffi(q)\frac{a(x)t}{\left(\pi|x-q|^2+a^2(x)t^2 \right)^{\frac{d+1}{2}}}dq,
 \end{align*}
 where $\Gamma$ is the Euler gamma-function. We see that the multiplicative perturbation $a(x)$ of the fractional Laplacian contributes actually to the time parameter in the definition of the family $(F(t))_{t\geq0}$. This motivates the result of the following subsection.
\end{example}

%%%%%%%%%%%%%%%%%%%%%%%%%%%%%%%%%%%%%%%%%%%%%%%%%%%%%%%%%%%%%%%%%%%%%%%%%%%%%%%%%%%%%%%%%%%%%%%%%%%%%%%%%%%%%%

\subsection{Chernoff approximations for  multiplicative perturbations of a generator}\label{subsec_mult}
Let $Q$ be a metric space. Consider the Banach space $X=C_b(Q)$ of bounded continuous functions on   $Q$ with supremum-norm  $\|\cdot\|_\infty$. Let   $(T_t)_{t \geq 0}$ be a strongly continuous semigroup on $X$  with generator  $(L,\Dom(L))$. Consider a  function  $a\in C_b(Q)$ such that $a(q)>0$ for all $q\in Q$. Then the space $X$ is invariant under the multiplication operator $a$, i.e.   $a(X)\sbs X$. Consider the operator 
$\Sharp{L}$, defined for all $\ffi\in\Dom(\Sharp{L})$ and all $q\in Q$  by 
\begin{equation}\label{1:2:eq:Lperturbed}
\Sharp{L}\varphi(q):=a(q)(L\varphi)(q),\quad \text{where}\quad\Dom(\Sharp{L}):=\Dom(L).
\end{equation}
\begin{assumption}\label{1:ass:multPert}
We assume that  $(\Sharp{L},\Dom(\Sharp{L}))$ generates a strong\-ly continuous semigroup (which is denoted by $(\Sharp{T}_t)_{t \geq 0}$) on the  Banach space $X$.
\end{assumption}
Some conditions assuring the existence and strong continuity of the semigroup $(\Sharp{T}_t)_{t \geq 0}$ can be found,  e.g., in~\cite{MR0201983,MR0338832}. The operator  $\Sharp{L}$ is called a \emph{multiplicative perturbation} of the generator  $L$ and the semigroup $(\Sharp{T}_t)_{t \geq 0}$, generated by $\Sharp{L}$, is called a \emph{semigroup with the multiplicatively perturbed } with the function  $a$ \emph{generator}.  The following result has been shown in~\cite{Habil} (cf. \cite{technomag,MR2999096}).
\begin{theorem}\label{1:thm:multiplicPert-Cb}
Let Assumption~\ref{1:ass:multPert} hold.
Let  $(F(t))_{t\ge0}$ be a strongly continuous family\footnote{The family $(F(t))_{t\ge0}$ of bounded linear operators on a Banach space $X$ is called strongly continuous if $\lim_{t\to t_0}\|F(t)\ffi-F(t_0)\ffi\|_X=0$ for all $t,t_0\geq0$ and all $\ffi\in X$.} of bounded linear operators on the Banach space  $X$, which is Chernoff equivalent to the semigroup  $({T}_t)_{t \ge 0}$.  Consider  the family  of operators $(\Sharp{F}(t))_{t\ge0}$ defined on $X$ by
\begin{equation}\label{1:eq:F(t) for multPert}
\Sharp{F}(t)\ffi(q):=(F(a(q)t)\ffi)(q)\qquad\text{ for all }\,\ffi\in X,\,\,\, q\in Q.
\end{equation}
The operators $\Sharp{F}(t)$ act on the space $X$,  the family  $(\Sharp{F}(t))_{t\ge0}$ is again strongly continuous and is Chernoff equivalent to the semigroup $(\Sharp{T}_t)_{t \ge 0}$  with multiplicatively perturbed with the function  $a$ generator, i.e.  the Chernoff approximation  
\begin{equation*}
\Sharp{T}_t\ffi=\lim_{n\to\infty}\big[\Sharp{F}(t/n)\big]^n\ffi
\end{equation*}
is valid for all $\ffi\in X$ locally uniformly with respect to  $t\geq0$. 
\end{theorem}
\begin{remark}\label{1:prop:multPert}
(i) The statement of Theorem~\ref{1:thm:multiplicPert-Cb}  remains true  for the following Banach spaces (cf.~\cite{Habil}):
%\begin{enumerate}

\noindent (a) $\quad X=C_\infty(Q):=\left\{\varphi\in C_b(Q)\,:\, \lim_{\rho(q,q_0)\to\infty}\varphi(q)=0 \right\},$ where $q_0$  is an arbitrary fixed point of   $Q$ and the metric space $Q$ is unbounded with respect to its metric $\rho$;

\noindent (b) $\quad X=C_0(Q):=\big\{ \ffi\in C_b(Q)\,:\, \forall\,\eps>0\,\,\exists\,\text{ a compact }\, K^\eps_\ffi\sbs Q\,\text{ such that }\, |\ffi(q)|<\eps\,$ $\text{ for all }\, q\notin K^\eps_\ffi    \big\}$, where the metric space $Q$ is assumed to be locally compact.%\footnote{If $Q=\cRd$, we have $C_\infty(Q)=C_0(Q)$, and we use the notation $C_\infty(\cRd)$ for this space. If $Q=(0,\infty)$, we have $C_0(Q)=C_\infty(Q)\cap\{  \ffi:\,\liml_{x\searrow0}\ffi(x)=0  \}$ and $C_\infty((0,\infty))=C_0([0,\infty))$. In general, it is not assumed in the definition of $C_\infty(Q)$ that $Q$ is locally compact.}.
%\end{enumerate}

\noindent(ii) As it follows from the proof of Theorem~\ref{1:thm:multiplicPert-Cb}, if  $\lim_{t\to0}\big\|
\frac{F(t)\ffi-\ffi}{t}-L\ffi\big\|_X=0$ for all  $\ffi\in D$ then also $\lim_{t\to0}\big\|
\frac{\Sharp{F}(t)\ffi-\ffi}{t}-\Sharp{L}\ffi\big\|_X=0$ for all  $\ffi\in D$.
\end{remark}

\begin{corollary}\label{1:cor:multPert}
Let   $(X_t)_{t\ge0}$ be a Markov process with the state space  $Q$  and transition probability $P(t,q,dy)$. Let the corresponding semigroup     $({T}_t)_{t\ge0}$,
$$
T_t\varphi(q)=\mathbb{E}^q\left[\varphi(X_t)\right]\equiv\intl_{Q}\varphi(y)P(t,q,dy),
$$
be strongly continuous on the Banach space  $X$, where $X=C_b(Q)$,  $X=C_\infty(Q)$ or $X=C_0(Q)$, and Assumption~\ref{1:ass:multPert} hold. Then by Theorem~\ref{1:thm:multiplicPert-Cb} and Remark~\ref{1:prop:multPert} the family   $(\Sharp{F}(t))_{t\ge0}$ defined by
$$
\Sharp{F}_t\varphi(q):=\intl_{Q}\varphi(y)P({a(q)t},q,dy),
$$
is strongly continuous and is Chernoff equivalent to the semigroup  $(\Sharp{T}_t)_{t \geq 0}$ with multiplicatively perturbed (with the function  $a$) generator. Therefore, the following Chernoff approximation %Lagrangian Feynman formula
 is true for all $t>0$ and all $q_0\in Q$:
\begin{equation} \label{1:eq:LFF-multPert}
\begin{aligned}
\Sharp{T}_t\varphi(q_0)=
\lim_{n\to\infty}\intl_{Q}\cdots \intl_{Q}\varphi(q_n)
P({a(q_0)t/n},&q_0,dq_1)P({a(q_1)t/n},q_1,dq_2)\times\cdots\\
&
\times P({a(q_{n-1})t/n},q_{n-1},dq_n),
\end{aligned}
\end{equation}
where the order of integration is from $q_n$ to $q_1$ and the convergence is uniform with respect to $q_0\in Q$ and locally uniform with respect to $t\geq0$.
\end{corollary}

\begin{remark}
A multiplicative perturbation of the generator of a Markov process is equivalent to some randome time change of the process  (see \cite{MR0100919}, \cite{MR0137154}, \cite{MR838085}). Note that   $\Sharp{P}(t,q,dy):=P({a(q)t},q,dy)$  is not a transition probability  any more.  Nevertheless, if the  transition probability  $P(t,q,dy)$ of the original process is known, formula~\eqref{1:eq:LFF-multPert} allows to approximate the unknown transition probability of the modified process. %Some explicilty known transitional densities can be found, e.g., in  \cite{MR1912205,MR2042661}.  
\end{remark}

%%%%%%%%%%%%%%%%%%%%%%%%%%%%%%%%%%%%%%%%%%%%%%%%%%%%%%%%%%%%%%%%%%%%%%%%%%%%%%%%%%%%%%%%%%%%%%%%%%%%%%%%%%%%%%

\subsection{Chernoff approximations for semigroups generated by processes in a domain with prescribed behaviour at the boundary of / outside the domain}\label{subsec_BC}
%\subsection{Chernoff approximation for semigroups generated by reflected diffusions}
Let $(\xi_t)_{t\geq0}$ be a Markov process in $\cRd$. Assume that the corresponding semigroup $(T_t)_{t\geq0}$ is strongly continuous on some Banach space $X$ of functions on $\cRd$, e.g. $X=C_\infty(\cRd)$ or $X=L^p(\cRd)$, $p\in[1,\infty)$. Let $(L,\Dom(L))$ be the generator of $(T_t)_{t\geq0}$ in $X$. Assume that a Chernoff approximation of $(T_t)_{t\geq0}$ via a family $(F(t))_{t\geq0}$ is already known (and hence we have a core $D$ for $L$ such that $\lim_{t\to0}\big\| \frac{F(t)\ffi-\ffi}{t}-L\ffi\big\|_X=0$ for all  $\ffi\in D$). Consider now a domain $\Omega\subset\cRd$. Let $(\xi_t)_{t\geq0}$ start in $\Omega$ and impose some reasonable ``Boundary Conditions'' (\textbf{BC}), i.e. conditions on the behaviour of $(\xi_t)_{t\geq0}$ at the boundary $\pd\Omega$, or (if the generator $L$ is non-local) outside $\Omega$. This procedure gives rise to a Markov process  in $\Omega$ which we denote by $(\xi^*_t)_{t\geq0}$. 
In some cases, the corresponding semigroup $(T^*_t)_{t\geq0}$ is strongly continuous on some Banach space $Y$ of functions on $\Omega$ (e.g. $Y=C(\overline{\Omega})$, $Y=C_0(\Omega)$ or $Y=L^p(\Omega)$, $p\in[1,\infty)$). 
%Let  $(L^*,\Dom(L^*))$ be the generator of $(T^*_t)_{t\geq0}$ in $Y$. 
The question arises: \emph{how to construct a Chernoff approximation of $(T^*_t)_{t\geq0}$ on the base of the family $(F(t))_{t\geq0}$, i.e. how to incorporate \textbf{BC} into a Chernoff approximation}? A possible strategy to answer this question is to construct a proper extension  $E^*$ of functions from $\Omega$ to $\cRd$ such that, first, $E^*\,: Y\to X$ is a linear contraction and, second,  there exists a core $D^*$ for the generator $(L^*,\Dom(L^*))$  of $(T^*_t)_{t\geq0}$ with $E^*(D^*)\subset D$. Then it is easy to see that the family $(F^*(t))_{t\geq0}$ with 
\begin{align}\label{BC-ChAp}
F^*(t):=R_t\circ F(t)\circ E^*
\end{align}
 is Chernoff equivalent to the semigroup $(T^*_t)_{t\geq0}$. Here $R_t$ is, in most cases, just the restriction of functions from $\cRd$ to $\Omega$, and, for the case of Dirichlet \textbf{BC}, it is a multiplication  with a proper cut-off function $\psi_t$ having support in $\Omega$ such that $\psi_t\to 1_{\Omega}$ as $t\to0$ (see \cite{MR3903609,MR2729591}). This stategy has been successfully realized in the following cases (note that extensions $E^*$ are obtained in a constructive way and can be implemented in numerical schemes):

%\begin{enumerate}
%\item[(1)]

\textbf{Case 1:} $X=C_\infty(\cRd)$,  $(\xi_t)_{t\geq0}$ is a Feller process whose generator $L$ is given  by~\eqref{2:eq:g35}  with  $A$, $B$, $C$ of the class $C^{2,\alpha}$, $A$ satisfies~\eqref{eq:CP:uniformEllipt}, and either $N\equiv0$ or $N\neq0$ and the non-local term of $L$ is a relatively bounded perturbation of the local part of $L$ with some extra assumption on jumps of the process (see details in~\cite{MR3903609,MR2729591}). The family $(F(t))_{t\geq0}$ is given by~\eqref{eq:F(t)-NonLocCP} (see also~\eqref{F(t)via2conv},  or~\eqref{F(t)-Gaus} in the corresponding particular cases) and $D=C^{2,\alpha}_c(\cRd)$.   Further, $\Omega$ is a bounded $C^{4,\alpha}-$smooth  domain,  $Y=C_0(\Omega)$, \textbf{BC} are the homogeneous Dirichlet boundary/external conditions corresponding to killing of the process upon leaving the domain $\Omega$. A proper extension $E^*$ has been constructed in~\cite{MR2860750}, and it maps $\Dom(L^*)\cap C^{2,\alpha}(\overline{\Omega})$ into $D$.

One can further simplify the  Chernoff approximation constructed via the family $(F^*(t))_{t\geq0}$ of~\eqref{BC-ChAp} and show that the following Feynman formula solves the considered Cauchy--Dirichlet problem (see~\cite{MR3903609}):
\begin{align}\label{eq:generalFFforCDP}
T^*_t\ffi(x_0)=\liml_{n\to\infty}\intl_{\Omega}\ldots\intl_{\Omega}\intl_{\Omega}\ffi(x_n)\,\nu^{x_{n-1}}_{t/n}(dx_n)\nu^{x_{n-2}}_{t/n}(dx_{n-1})\cdots \nu^{x_0}_{t/n}(dx_1).
\end{align}
The convergence in this formula is however only locally uniform with respect to $x_0\in \Omega$ (and locally uniform with respect to  $\geq0$).
Similar results hold also for non-degenerate diffusions in domains of a compact Riemannian manifold $M$ with homogeneous Dirichlet \textbf{BC} (see, e.g.~ \cite{MR2423533}), what can be shown by combining approaches described in Subsections~\ref{subsec_additive},~\ref{subsec_mult},~\ref{subsec_BC} and using families $(F(t))_{t\geq0}$ of~\cite{MR2276523} which are Chernoff equivalent to the heat semigroup on $C(M)$.

%\bigskip

%\item[(2)]

\textbf{Case 2:}  $X=C_\infty(\cRd)$, $(\xi_t)_{t\geq0}$ is a Brownian motion, the family $(F(t))_{t\geq0}$ is the heat semigroup~\eqref{heatSG}  (hence $D=\Dom(L)$), $\Omega$ is a bounded  $C^\infty$-smooth domain, $Y=C(\overline{\Omega})$, \textbf{BC} are the Robin boundary conditions
\begin{align}\label{RobinBC}
\frac{\pd\ffi}{\pd\nu}+\beta\ffi=0\quad\text{ on }\,\,\pd\Omega,
\end{align}
where $\nu$ is the outer unit normal, $\beta$ is a smooth bounded nonnegative function on $\pd\Omega$. A proper extension $E^*$ (and the corresponding Chernoff approximation itself) has been constructed in~\cite{MR2541275}, and it maps $D^*:=\Dom(L^*)\cap C^\infty(\overline{\Omega})$ into the $\Dom(L)$. 

This result can be further generalized for the case of  diffusions, using the techniques of Subsections~\ref{subsec_additive} and~\ref{subsec_mult}. This will be demonstrated in Example~\ref{Example-Robin}. Note, however, that the extension $E^*$ of~\cite{MR2541275} maps $D^*$ into the set of functions which do not belong to $C^2(\cRd)$. Hence it is not possible to use the family~\eqref{F(t)-Gaus} (and $D=C^{2,\alpha}_c(\cRd)$) in a straightforward manner for approximation of diffusions with Robin \textbf{BC}.  
%\end{enumerate}

\begin{example}\label{Example-Robin}
Let $X=C_\infty(\cRd)$. Consider $(L_1,\Dom(L_1))$ being the closure of $\left(\frac12\Delta, S(\cRd)\right)$ in $X$.
Then $\Dom(L_1)$ is continuously embedded in $C^{1,\alpha}(\cRd)$ for every $\alpha\in(0,1)$ by Theorem~3.1.7 and Corollary~3.1.9~(iii) of~\cite{MR3012216}, and $(L_1,\Dom(L_1))$ generates a $C_0$-semigroup  $(T_1(t))_{t\geq}$ on $X$, this is the heat semigroup given by~\eqref{heatSG}. Let $a\in C_b(\cRd)$  be such that $a(x)\geq a_0$ for some $a_0>0$ and all $x\in\cRd$. Then $(\Sharp{L_1},\Dom(L_1))$, $\Sharp{L_1}\ffi(x):=a(x)L_1\ffi(x)$, generates a $C_0$-semigroup $(\Sharp{T_1}(t))_{t\geq0}$ on $X$ by~\cite{MR0201983}. Therefore, the family $(\Sharp{F}(t))_{t\geq0}$ with
\begin{align*}
\Sharp{F}(t)\ffi(x):=\int_\cRd \ffi(y) P(a(x)t,x,dy),
\end{align*}
where $P(t,x,dy)$ is given by~\eqref{heat-P}, is Chernoff equivalent to $(\Sharp{T_1}(t))_{t\geq0}$  by Corollary~\ref{1:cor:multPert}. And $\left\|\frac{\Sharp{F}(t)\ffi-\ffi}{t}-\Sharp{L_1}\ffi\right\|_X\to0$ as $t\to0$ for each $\ffi\in\Dom(L_1)$. Let now $C\in C_b(\cRd)$ and $B\in C_b(\cRd;\cRd)$. Then the operator $(L,\Dom(L))$,
\begin{align*}
L\ffi(x):=\frac{a(x)}{2}\Delta\ffi(x)-B(x)\cdot\nabla\ffi(x)-C(x)\ffi(x),\quad\Dom(L):=\Dom(L_1),
\end{align*}
(obtained by a relatively bounded additive perturbation of $(\Sharp{L_1},\Dom(L_1))$) generates a $C_0$-semigroup $(T_t)_{t\geq0}$ on $X$ (e.g., by Theorem 4.4.3 of~\cite{MR1873235}). Motivated by Subsection~\ref{subsec_mult} and the view of the translation semigroup, consider the family $(F_2(t))_{t\geq0}$ of contractions on $X$ given by $F_2(t)\ffi(x):=\ffi(x-tB(x))$. Then, for all $\ffi\in\Dom(L_1)\subset C^{1,\alpha}(\cRd)$, holds %(with some constant $C$) 
$ %\begin{align*}
\left\|\frac{F_2(t)\ffi-\ffi}{t}+B\cdot\nabla\ffi\right\|_X\leq \text{const}\cdot t^\alpha|B|^{\alpha+1}\to0,\quad t\to0,
$  and $\|F_2(t)\|\leq1$ for all $t\geq0$. 
%\end{align*}
Therefore, by Theorem~\ref{1:cor:additivePert}, the family $(F(t))_{t\geq0}$ with 
\begin{align*}
F(t)\ffi(x):&=\left[e^{-tC}\circ F_2(t)\circ \Sharp{F}(t)\right]\ffi(x)\\
&\equiv e^{-tC(x)}\int_{\cRd}\ffi(y) P(a(x-tB(x))t,x-tB(x),dy)
\end{align*}
is Chernoff equivalent to the semigroup $(T_t)_{t\geq0}$.
Let now $\Omega$ be a bounded  $C^\infty$-smooth domain, $Y=C(\overline{\Omega})$. Consider $(L^*,\Dom(L^*))$
in $Y$ with\footnote{Here we consider the Robin \textbf{BC}~\eqref{RobinBC} given in a weaker form via the first Green's formula;  $d\sigma$ is the surface measure on $\pd\Omega$.}
\begin{align*}
&\Dom(L^*):=\bigg\{\ffi\in Y\cap H^1(\Omega)\,\,:\,\, L\ffi\in Y, \\
&
\qquad\qquad\qquad\int_\Omega\Delta\ffi u dx+     \int_\Omega\nabla\ffi\nabla udx+\int_{\pd\Omega}\beta\ffi u d\sigma=0\,\,\forall\,u\in H^1(\Omega)\bigg\},\\
&
L^*\ffi:=L\ffi,\quad\forall\,\ffi\in\Dom(L^*).
\end{align*}
 Then $(L^*,\Dom(L^*))$ generates a $C_0$-semigroup $(T^*_t)_{t\geq0}$ on $Y$ (cf.~\cite{MR2812574}). Consider $R\,:\,X\to Y$ being the restriction of a function from $\cRd$ to $\overline{\Omega}$. Consider the extension $E^*\,:\,Y\to X$ constructed in~\cite{MR2541275}. This extension  is a linear contraction, obtained via an orthogonal reflection at the boundary and multiplication with a suitable  cut-off function, whose behaviour at $\pd\Omega$ is prescribed (depending on $\beta$) in such a way that the weak Laplacian of the extension $E^*(\ffi)$ is continuous for each $\ffi\in D^*:=\Dom(L^*)\cap C^\infty(\Omega)$ and $E^*(D^*)\subset \Dom(L_1)$. We omit the explicit description of $E^*$, in order to avoid corresponding technicalities. We consider the family  $(F^*(t))_{t\geq0}$ on $Y$ given by $F^*(t):=R\circ F(t)\circ E^*$, i.e.
 $$
 F^*(t)\ffi(x):=e^{-tC(x)}\int_{\cRd}E^*[\ffi](y) P(a(x-tB(x))t,x-tB(x),dy),\qquad x\in\overline{\Omega}.
 $$ 
 Then $F^*(0)=\Id$,  $\|F^*(t)\|\leq e^{t\|C\|_\infty}$, and we have for all $\ffi\in D^*$ 
\begin{align*}
\liml_{t\to0}\left\|\frac{F^*(t)\ffi-\ffi}{t}-L^*\ffi  \right\|_Y&=\liml_{t\to0}\left\|R\circ\left(\frac{F(t)E^*[\ffi]-E^*[\ffi]}{t}-L E^*[\ffi]\right)  \right\|_Y\\
&
\leq\liml_{t\to0}\left\|\frac{F(t)E^*[\ffi]-E^*[\ffi]}{t}-L E^*[\ffi]  \right\|_X =0.
\end{align*}
Therefore, the family $(F^*(t))_{t\geq0}$ is Chernoff equivalent to the semigroup $(T^*_t)_{t\geq0}$ by Theorem~\ref{ChernoffTheorem}, i.e. 
% \begin{equation}\label{0:eq:ChernoffFormula}
$T^*_t\ffi=\liml_{n\to\infty}[F^*(t/n)]^n\ffi$ for each $\ffi\in Y$ locally uniformly with respect to $t\geq0$.
%\end{equation}

 \end{example}

\subsection{Chernoff approximations for subordinate  semigroups }\label{subsec_subordinate}
One of the ways to construct strongly continuous semigroups is given by the procedure of subordination. From  two ingredients: an original $C_0$ contraction semigroup $(T_t)_{t\geq0}$ on a Banach space $X$ and a convolution semigroup\footnote{A family $(\eta_t)_{t\geq0}$ of bounded Borel measures on $\cR$ is called a convolution semigroup   if $\eta_t(\cR)\leq 1$ for all $t\geq0$,  $\eta_t*\eta_s=\eta_{t+s}$ for all $t,s\geq0$, $\eta_0=\delta_0$, and $\eta_t\to\delta_0$ vaguely as $t\to0$, i.e. $\lim_{t\to0}\int_{\cR}\ffi(x)\eta_t(dx)=\int_{\cR}\ffi(x)\delta_0(dx)\equiv \ffi(0)$ for all $\ffi\in C_c(\cR)$. A convolution semigroup $(\eta_t)_{t\geq0}$ is supported by $[0,\infty)$ if $\supp\eta_t\subset[0,\infty)$ for all $t\geq0$.}  $(\eta_t)_{t\geq0}$ supported by $[0,\infty)$, this procedure produces  the $C_0$ contraction semigroup $(T^f_t)_{t\geq0}$ on $X$  with  
$$
 T^f_t\ffi:=\int_0^\infty T_s\ffi\,\eta_t(ds),\quad\forall\,\,\ffi\in X.
$$
If the semigroup $(T_t)_{t\geq0}$ corresponds to a stochastic process $(X_t)_{t\geq0}$, then  subordination is a random time-change of $(X_t)_{t\geq0}$ by an independent increasing L\'{e}vy process (subordinator) with distributions $(\eta_t)_{t\geq0}$. If $(T_t)_{t\geq0}$ and  $(\eta_t)_{t\geq0}$ both are known explicitly, so is $(T^f_t)_{t\geq0}$. But if, e.g., $(T_t)_{t\geq0}$ is not known, neither $(T^f_t)_{t\geq0}$ itself, nor even the generator of $(T^f_t)_{t\geq0}$ are known explicitly any more. This impedes the construction of a family $(F(t))_{t\geq0}$  with a prescribed  (but unknown explicitly) derivative at $t=0$. This difficulty is overwhelmed below by construction of families $(\mathcal{F}(t))_{t\geq0}$ and $(\mathcal{F}_\mu(t))_{t\geq0}$  which incorporate approximations of the generator of $(T^f_t)_{t\geq0}$   itself. Recall that each convolution semigroup $(\eta_t)_{t\geq0}$ supported by $[0,\infty)$ corresponds to a Bernstein function $f$ via the Laplace transform  $\mathcal{L}$:  $\mathcal{L}[\eta_t]=e^{-tf}$ for all   $t>0$. Each Bernstein function $f$ is uniquely defined by a triplet $(\sigma,\lambda,\mu)$ with constants $\sigma,\lambda\geq0$ and  a Radon measure $\mu$  on $(0,\infty)$, such that $\int_{0+}^\infty\frac{s}{1+s}\mu(ds)<\infty$,
 through the representation
% \begin{equation}\label{eq:BernsteinFunc}
$f(z)=\sigma+\lambda z+\intl_{0+}^\infty (1-e^{-sz})\mu(ds),\quad\quad\any\,z\,: \RE z\geq 0.$
%\end{equation}
%Note that   $\eta_t(\cR)=1$ for all $ t\geq0$ if and only if $\sigma=0$ (i.e. there is no "killing", cf. \cite{MR3156646}).   
Let $(L,\Dom(L))$ be the generator of $(T_t)_{t\geq0}$ and $(L^f,\Dom(L^f))$ be the generator of $(T^f_t)_{t\geq0}$. Then each core for $L$ is also a core for $L^f$ and,   for $\ffi\in\D(L)$, the operator $L^f$ has the representation
\begin{equation*}\label{eq:L^f}
L^f\ffi=-\sigma\ffi+\lambda L\ffi+\intl_{0+}^\infty (T_s\ffi-\ffi)\mu(ds).
\end{equation*} 
Let  $(F(t))_{t\geq 0}$ be a family of contractions on $(X,\|\cdot\|_X)$ which is Chernoff equivalent to $(T_t)_{t\geq 0}$, i.e. $F(0)=\id$,  $\|F(t)\|\leq 1$ for all  $t\geq 0$ and there is a set $D\subset \D(L)$, which is a core for $L$, such that  $\lim_{t\to0}\big\|\frac{F(t)\ffi-\ffi}{t}-L\ffi\big\|_X=0$ for each $\ffi\in D$. The first candidate
for being Chernoff equivalent to $(T^f_t)_{t\geq0}$ could be  the family of operators $(F^*(t))_{t\geq 0}$ given by  $F^*(t)\ffi:=\int_0^\infty F(s)\ffi\,\eta_t(ds)$ for all $\ffi\in X$. However, its derivative at zero does not coincide with $L^f$ on $D$. Nevertheless, with suitable modification of $(F^*(t))_{t\geq 0}$,  Theorem~\ref{1:cor:additivePert} and the discussion below Theorem~\ref{ChernoffTheorem}, the following has been proved in~\cite{MR3804267}.
\begin{theorem}\label{thm:ChFam_knownSub-r}
  Let $m:(0,\infty)\to\mathbb{N}_0$ be a monotone function\footnote{ One can take, e.g., $m(t):=\left\lfloor 1/t  \right\rfloor=$  the largest integer $n\leq 1/t$. Recall that $\mathbb{N}_0:=\mathbb{N}\cup\{0\}.$ } such that 
   $m(t)\to+\infty$ as $t\to0$.
Let the mapping  $[F(\cdot/m(t))]^{m(t)}\ffi\,:\,[0,\infty)\to X $ be Bochner  measurable  as the mapping from  $([0,\infty),\mathcal{B}([0,\infty)),\eta^0_t)$ to $(X,\mathcal{B}(X))$  for each $t>0$ and each $\ffi\in X$.
% Let $f$ be a Bernstein function  given by a triplet $(\sigma,\lambda,\mu)$ through  the representation \eqref{eq:BernsteinFunc}  with associated convolution semigroup $(\eta_t)_{t\geq 0}$ supported by $[0,\infty)$.   

\textbf{Case 1:} 
Let $(\eta^0_t)_{t\geq 0}$ be the  convolution semigroup (supported by $[0,\infty)$) associated to the Bernstein function $f_0$ defined by the triplet $(0,0,\mu)$.  Assume that the corresponding operator semigroup $(S_t)_{t\geq0}$, $S_t\ffi:=\ffi*\eta^0_t$, is strong Feller\footnote{The  semigroup $(S_t)_{t\geq0}$ is   strong Feller 
 iff all the measures $\eta^0_t$ admit densities of the class $L^1([0,\infty))$ with respect to the Lebesgue measure (cf. Examle~4.8.21 of \cite{MR1873235}).}.
       Consider the family   $(\mathcal{F}(t))_{t\geq 0}$  of operators on $(X,\|\cdot\|_X)$ defined by $\mathcal{F}(0):=\id$ and
\begin{equation}\label{eq:mathcal(F)}
\mathcal{F}(t)\ffi:=  e^{-\sigma t}\circ   F(\lambda t)\circ \mathcal{F}_0(t)\ffi,\quad t>0,\, \ffi\in X,
\end{equation}
with $\mathcal{F}_0(0)=\id$ and\footnote{For any bounded operator $B$, its zero degree $B^0$ is considered to be the identity operator. For each $t>0$, a non-negative integer $m(t)$ and  a  bounded Bochner measurable  mapping   $[F(\cdot/m(t))]^{m(t)}\ffi\,:\,[0,\infty)\to X$, the integral in the right hand side of formula~\eqref{eq:F_0} is well defined.}
\begin{equation*}\label{eq:F_0}
\mathcal{F}_0(t)\ffi:=   \intl_{0+}^\infty \left[F( s/m(t))\right]^{m(t)}\ffi\,\eta^0_t(ds),\quad t>0,\, \ffi\in X.
\end{equation*}
The family $(\mathcal{F}(t))_{t\geq 0}$ is Chernoff equivalent to the semigroup $(T^f_t)_{t\geq 0}$, and hence
%%, the Feynman formula
$$
T^f_t \ffi=\lim_{n\to\infty} \big[\mathcal{F}(t/n)   \big]^n \ffi
$$
for all $\ffi\in X$ locally uniformly with respect to $t\geq  0$.

\textbf{Case 2:}   
Assume that the measure $\mu$ is bounded. 
 Consider a family   $(\mathcal{F}_\mu(t))_{t\geq 0}$  of operators on $(X,\|\cdot\|_X)$ defined for all $\ffi\in X$ and all $t\geq0$ by 
\begin{equation*}\label{eq:mathcal(F)mu}
\mathcal{F}_\mu(t)\ffi:=  e^{-\sigma t} F(\lambda t)\left( \ffi+ t \intl_{0+}^\infty (F^{m(t)}(s/m(t))\ffi-\ffi)\mu(ds)   \right).
\end{equation*}
The family $(\mathcal{F}_\mu(t))_{t\geq 0}$ is Chernoff equivalent to the semigroup $(T^f_t)_{t\geq 0}$, and hence
%%, the Feynman formula
$$
T^f_t \ffi=\lim_{n\to\infty} \big[\mathcal{F}_\mu(t/n)   \big]^n \ffi
$$
for all $\ffi\in X$ locally uniformly with respect to $t\geq  0$. 
\end{theorem}
The constructed families $(\mathcal{F}(t))_{t\geq0}$ and $(\mathcal{F}_\mu(t))_{t\geq0}$ can be used (in combination with the techniques of Sec.~\ref{subsec_additive}, Sec.~\ref{subsec_mult}, Sec.~\ref{subsec_BC} and results of~\cite{MR2927703,MR2276523}), e.g., to approximate semigroups generated by subordinate Feller diffusions on star graphs and Riemannian manifolds.
Note that the family~\eqref{eq:mathcal(F)} can be used when   the  convolution semigroup $(\eta^0_t)_{t\geq 0}$ is known explicitly.   This is the case of inverse Gaussian (including $1/2$-stable) subordinator,   Gamma subordinator and some others (see, e.g., \cite{MR1912205,MR3231629,MR2042661} for examples).

%\subsection{Applications to non-markovian evolution}

\subsection{Approximation of solutions of time-fractional evolution equations}\label{subsecFrac}
We are interested now in distributed order  time-fractional evolution equations of the form
\begin{align}\label{Appl:3:eq:fracFPK}
\mathcal{D}^\mu f(t)=Lf(t),
\end{align}
 where $(L,\Dom(L))$ is the generator of a $C_0-$ contraction semigroup $(T_t)_{t\geq0}$ on some Banach space $(X,\|\cdot\|_X)$  and $\mathcal{D}^\mu$ is the distributed order fractional derivative with respect to the time variable $t$:  
\begin{align*}
\mathcal{D}^\mu u(t):=\int_0^1 \frac{\pd^\beta}{\pd t^\beta} u(t)\mu(d\beta),\quad\text{ where } \quad \frac{\pd^\beta}{\pd t^\beta} u(t):=\frac{1}{\Gamma(1-\beta)}\int_0^t \frac{u'(r)}{(t-r)^\beta}dr,
\end{align*}
where $\mu$ is a finite Borel measure with $\supp\mu\in(0,1)$. Equations of such type are called \emph{time-fractional Fokker--Planck--Kolmogorov equations} (tfFPK-equations) and arise in the framework of continuous time random walks  and fractional kinetic theory (\cite{MR0260036,MR1809268,MR757002,MR3351502,MR1937584}).
As it is shown in papers \cite{MR2886388,MR2782245,MR3168478},  such tfFPK-equations are governing equations for stochastic processes which are  time-changed Markov processes, where the time-change $(E^\mu_t)_{t\geq0}$ arises as the first hitting time of level $t>0$ (or, equivalently, as the inverse process) %, i.e. $E^\mu_t:=\inf\left\{\tau\geq0\,:\,D^\mu_\tau>t \right\}$)
  for a  mixture  $(D^\mu_t)_{t\geq0}$ of independent stable subordinators with the  mixing measure $\mu$\footnote{Hence $(D^\mu_t)_{t\geq0}$ is a subordinator corresponding to the Bernstein function  $f^\mu(s):=\int_0^1 s^\beta\mu(d\beta)$,  $s>0$, and $E^\mu_t:=\inf\left\{\tau\geq0\,:\,D^\mu_\tau>t \right\}$.}.  Existence and uniqueness of solutions of initial and initial-boundary value problems for such tfFPK-equations are considered, e.g., in~\cite{MR3721899,MR3753615}.  The process $(E^\mu_t)_{t\geq0}$  is sometimes called \emph{inverse subordinator}. However,  note that it is not a Markov process. Nevertheless, $(E^\mu_t)_{t\geq0}$ posesses a nice  marginal density function $p^\mu(t,x)$   (with respect to the Lebesgue measure $dx$). It has been shown in~\cite{MR2886388,MR3168478} that the family of linear operators $(\mathcal{T}_t)_{t\geq0}$ from $X$ into $X$ given by
\begin{align}\label{Appl:3:eq:Tt-subordinated}
\mathcal{T}_t\ffi:=\int_0^\infty T_\tau\ffi\, p^\mu(t,\tau)\,d\tau,\quad\forall\,\ffi\in X,
\end{align}
is uniformly bounded,  strongly continuous,  and the function $f(t):=\mathcal{T}_t f_0$ is a solution of the Cauchy problem
\begin{align}\label{Appl:3:eq:CP}
&\mathcal{D}^\mu f(t)=Lf(t),\quad t>0,\nonumber\\
&
f(0)=f_0.
\end{align}
%Note that the family $(\mathcal{T}_t)_{t\geq0}$ does not posess the semigroup property. 
 This result shows that solutions of tfFPK-equations are a kind of subordination of  solutions of the corresponding time-non-fractional evolution equations with  respect to ``subordinators'' $(E^\mu_t)_{t\geq0}$. 
 %
 %And  respectively,  if a time-non-fractional evolution equation is a governing equation for a  Markov process then the related time-fractional evolution equation is a governing equation for a (already non-markovian) process which is a ``subordination'', i.e.  a time-change of this Markov process by means of $(E^\mu_t)_{t\geq0}$.
% 
  The non-Markovity of $(E^\mu_t)_{t\geq0}$ implies that the family $(\mathcal{T}_t)_{t\geq0}$ is not a semigroup any more.   
Hence we have no chances to construct Chernoff approximations for  $(\mathcal{T}_t)_{t\geq0}$.   Nevertheless, the following is true (see~\cite{MR3903609}).

\begin{theorem}\label{Appl:3:thm}
%Let $(X,\|\cdot\|_X)$ be a Banach space. Let $(L,\Dom(L))$ be the generator of a  strongly continuous contraction semigroup $(T_t)_{t\geq0}$ on $X$.
% Let $f_0\in\Dom(L)$.
Let the family $(F(t))_{t\geq0}$ of contractions on $X$ be 
 Chernoff equivalent to $(T_t)_{t\geq0}$. Let $f_0\in\Dom(L)$. Let the mapping $F(\cdot)f_0\,:\,[0,\infty)\to X$ be Bochner measurable as a mapping from $([0,\infty),\mathcal{B}([0,\infty)),dx)$ to $(X,\mathcal{B}(X))$.
   Let $\mu$ be a finite Borel measure with $\supp\mu\in(0,1)$ and the family $(\mathcal{T}_t)_{t\geq0}$ be given by formula~\eqref{Appl:3:eq:Tt-subordinated}. Let $f\,:\,[0,\infty)\to X$ be defined via $f(t):=\mathcal{T}_t f_0$.  For each $n\in\Nat$  define  the mappings  $f_n\,:\,[0,\infty)\to X$  by
\begin{align}\label{eq:fn}
f_n(t):=\int_0^\infty F^n(\tau/n)f_0\,p^\mu(t,\tau)\,d\tau.
\end{align}
 Then it holds  locally uniformly with respect to  $t\geq0$ that
$$
\|f_n(t)-f(t)\|_X\to0,\quad n\to\infty.
$$
\end{theorem}
Of course, similar approximations are valid also in the case of ``ordinary subordination'' (by a L\'{e}vy subordinator) considered in Sec.~\ref{subsec_subordinate}. 
Note also that there exist different Feynman-Kac formulae for the Cauchy problem~\eqref{Appl:3:eq:CP}. In particular, the function
\begin{align}\label{Appl:3:eq:f-FKF}
f(t,x):=\E\left[f_0\left(\xi\left(E^\mu_t\right)\right)\,\,|\,\,\xi(E^\mu_0)=x\right],
\end{align}
where $(\xi_t)_{t\geq0}$ is a Markov process with generator $L$, solves the Cauchy problem~\eqref{Appl:3:eq:CP} (cf. Theorem~3.6 in \cite{MR2886388}, see also \cite{MR3753615}). Furthermore, the considered   equations (with $\mu=\delta_{\beta_0}$, $\beta_0\in(0,1)$) are related to some time-non-fractional evolution equations of higher order (see, e.g., \cite{MR2491905,MR2965747}). Therefore, the approximations  $f_n$ constructed in Theorem~\ref{Appl:3:thm}   can be used simultaneousely to  approximate  path integrals appearing in different sto\-chas\-tic representations of the same function $f(t,x)$ and to  approximate  solutions of corresponding time-non-fractional evolution equations of higher order.

\begin{example}\label{example-FracFeller}
Let $\mu=\delta_{1/2}$, i.e. $\mathcal{D}^\mu $ is the Caputo derivative of $1/2$-th order and $(E^{1/2}_t)_{t\geq0}$ is a $1/2$-stable inverse subordinator whose marginal probability density is known explicitly: %(see Cor.~3.1 in \cite{MR2074812} and the discussion after  Lemma~3 in \cite{MR2726092})
%\begin{align*}
$p^{1/2}(t,\tau)=\frac{1}{\sqrt{\pi t}}e^{-\frac{\tau^2}{4t}}$.
%\end{align*}
 Let $X=C_\infty(\cRd)$ and $(L,\Dom(L))$ be the Feller generator given by~\eqref{2:eq:g35}. Let all the assumptions of Theorem~\ref{2:thm:HFF}  be fulfilled. Hence we can use the family $(F(t))_{t\geq0}$ given by~\eqref{eq:F(t)-NonLocCP} (or by~\eqref{F(t)-Gaus} if $N\equiv0$). Therefore, by Theorem~\ref{2:thm:HFF} and Theorem~\ref{Appl:3:thm}, the following Feynman formula solves the Cauchy problem~\eqref{Appl:3:eq:CP}:
 \begin{align*}
f(t,x_0)=\liml_{n\to\infty}\intl_0^\infty\intl_{\cRd}\ldots\intl_{\cRd}\frac{1}{\sqrt{\pi t}}e^{-\frac{\tau^2}{4t}}\ffi(x_n)\,\nu^{x_{n-1}}_{\tau/n}(dx_n)\cdots \nu^{x_0}_{\tau/n}(dx_1)d\tau.
\end{align*} 
\end{example}

\subsection{Chernoff approximations for Schr\"{o}dinger groups}\label{subsec_Schr}
\textbf{Case 1: PDOs.}
In Sec.~\ref{subsec_Feller}, we have used the technique of pseudo-differential operators (PDOs). Namely, (with a slight modification of notations) we have considered  operator semigroups $(e^{-t\widehat{H}})_{t\geq0}$ generated by  PDOs $-\widehat{H}$ with symbols $-H$ (see formula~\eqref{PDO_L}). We have approximated semigroups via families  of PDOs $(F(t))_{t\geq0}$ with symbols $e^{-tH}$, i.e. $F(t)=\widehat{e^{-tH}}$. Note again that $e^{-t\widehat{H}}\neq \widehat{e^{-tH}}$ in general. It was established in Theorem~\ref{2:thm:HFF} that
\begin{align}\label{HFF-SG}
e^{-t\widehat{H}}=\liml_{n\to\infty}\left[ \widehat{e^{-tH/n}}\right]^n
\end{align}
for a class of symbols $H$ given by~\eqref{2:eq:g116}. The same approach can be used to construct Chernoff approximations for Schr\"{o}dinger groups  $(e^{-it\widehat{H}})_{t\in\cR}$ describing quantum evolution of systems obtained by a quantization of classical systems with Hamilton functions $H$. Namely, it holds under certain conditions
\begin{align}\label{HFF-Schr}
e^{-it\widehat{H}}=\liml_{n\to\infty}\left[ \widehat{e^{-itH/n}}\right]^n.
\end{align}
On a heuristic level, such approximations have been considered already in works~\cite{MR0479157,MR631347}. A rigorous mathematical treatment and some conditions, when~\eqref{HFF-Schr} holds, can be found in~\cite{MR1927359}. Note that right hand sides of both~\eqref{HFF-SG} and~\eqref{HFF-Schr} can be interpreted as phase space Feynman path integrals~\cite{MR0479157,MR631347,MR2863557,MR3455669,MR1927359}.

\textbf{Case 2: ``rotation''.}
Another approach to construct Chernoff approximations for Schr\"{o}dinger groups $(e^{itL})_{t\in\cR}$  is based on a kind of ``rotation'' of families $(F(t))_{t\geq0}$ which are Chernoff equivalent to semigroups $(e^{tL})_{t\geq0}$ (see~\cite{MR3490776}). Namely, let $(L,\Dom(L))$ be a self-adjoint operator in a Hilbert space $X$ which generates a $C_0-$semigroup $(e^{tL})_{t\geq0}$ on $X$. Let a family $(F(t))_{t\geq0}$ be Chernoff equivalent\footnote{Actually, $(F(t))_{t\geq0}$ does not need to fulfill the condition (ii) of the Chernoff Theorem~\ref{ChernoffTheorem} in this construction. } to $(e^{tL})_{t\geq0}$. Let the operators $F(t)$ be self-adjoint for all $t\geq0$. Then the family $(F^*(t))_{t\geq0}$,
\begin{align*}
F^*(t):=e^{i(F(t)-\Id)},
\end{align*}
is Chernoff equivalent to the Schr\"{o}dinger (semi)group $(e^{itL})_{t\geq0}$. Indeed, $F^*(0)=\Id$, $\|F^*(t)\|\leq1$ since all $F^*(t)$ are unitary operators, and $(F^*)'(0)=iF'(0)$.  Hence the following Chernoff approximation holds
\begin{align}\label{ChAp-Vania}
e^{itL}\ffi=\liml_{n\to\infty}e^{in(F(t/n)-\Id)}\ffi,\qquad\forall\,\,\ffi\in X.
\end{align}
Since all $F(t)$ are bounded operators, one can calculate $e^{in(F(t/n)-\Id)}$ via Taylor expansion or via formula~\eqref{formula-E}. Let us illustrate this approach with the following example.
\begin{example}\label{example-rotation}
Consider the function $H$ given by~\eqref{2:eq:g116}, Assume that $H$ does not depend on $x$, i.e. $H=H(p)$, and $H$ is real-valued (hence $B\equiv0$ and $N(dy)$ is symmetric). Such symbols $H$ correspond to symmetric L\'{e}vy processes. It is well-known\footnote{See, e.g., Example~4.7.28 in \cite{MR1873235}.} that the closure $(L,\Dom(L))$ of $(-\widehat{H},C^\infty_c(\cRd))$ generates a $C_0-$ semigroup $(T_t)_{t\geq0}$ on $L^2(\cRd)$; operators $T_t$ are self-adjoint and  coincide with operators $F(t)$ given in~\eqref{F(t)-PDO}, i.e. $T_t=\widehat{e^{-tH}}$ on $C^\infty_c(\cRd)$. Therefore, the Chernoff approximation~\eqref{ChAp-Vania} holds for the Schr\"{o}dinger (semi)group $(e^{itL})_{t\geq0}$ resolving  the Cauchy problem
$$
-i\frac{\pd f}{\pd t}(t,x)=Lf(t,x),\qquad f(0,x)=f_0(x)
$$
in $L^2(\cRd)$ with $L$ being the generator of a symmetric L\'{e}vy process. Note that this class of generators contains differential operators with constant coefficients (with $H(p)=C+iB\cdot p+p\cdot Ap$), fractional Laplacians (with $H(p):=|p|^\alpha$, $\alpha\in(0,2)$) and relativistic Hamiltonians (with $H(p):=\sqrt[\alpha]{|p|^\alpha+m}$, $m>0$, $\alpha\in(0,2]$). Assume additionally that $H\in C^\infty(\cRd)$. Then $F(t)\,:\, S(\cRd)\to S(\cRd)$ and it holds on $S(\cRd)$ (with $\mathcal{F}$ and $\mathcal{F}^{-1}$ being Fourier and inverse Fourier transforms respectively):
\begin{align*}
&F(t)=\mathcal{F}^{-1}\circ e^{-tH}\circ\mathcal{F};\qquad\qquad\qquad in(F(t/n)-\Id)= \mathcal{F}^{-1}\circ\left( in\left(e^{-tH/n}-1\right)\right)\circ\mathcal{F};\\
&
\left[F^*(t/n) \right]^n=  e^{in(F(t/n)-\Id)}=\suml_{k=0}^\infty\frac{1}{k!}\left(\mathcal{F}^{-1}\circ\left( in\left(e^{-tH/n}-1\right)\right)\circ\mathcal{F}\right)^k\\
&
\qquad\qquad=\mathcal{F}^{-1}\circ\left[\suml_{k=0}^\infty\frac{1}{k!}\left( in\left(e^{-tH/n}-1\right)^k\right)\right]\circ\mathcal{F}=\mathcal{F}^{-1}\circ \exp\left\{in\left(e^{-tH/n}-1\right)\right\}\circ\mathcal{F}.
\end{align*}
Therefore, $\left[F^*(t/n) \right]^n$ is a PDO with symbol $\exp\left\{in\left(e^{-tH/n}-1\right)\right\}$ on $S(\cRd)$. Hence we have obtained the following representation for the Schr\"{o}dinger (semi)group $(e^{itL})_{t\geq0}$:
\begin{align}\label{eq:Schr-Appr}
&e^{itL}\ffi(x)=\liml_{n\to\infty}(2\pi)^{-d}\intl_{\cRd}\intl_{\cRd}e^{ip\cdot(x-q)}\exp\left\{in\left(e^{-tH(p)/n}-1\right)\right\}\ffi(q)\,dqdp,
\end{align}
for all $\ffi\in S(\cRd)$ and all $x\in\cRd$. The convergence in~\eqref{eq:Schr-Appr} is in $L^2(\cRd)$ and is locally uniform with respect to $t\geq0$. 
\end{example}

\textbf{Case 3: shifts and averaging.} One more approach to construct Chernoff approximations for semigroups and Schr\"{o}dinger groups generated by differential and pseudo-differential operators is based on shift operators (see~\cite{MR3588869,MR3871530}), averaging (see~\cite{MR3819810,MR3479996}) and their combination (see~\cite{MR3819810,MR3844136}). Let us demonstrate this method by means of  simplest examples. So, consider $X=C_\infty(\cR)$ or $X=L^p(\cR)$, $p\in[1,\infty)$. Consider $(L,\Dom(L))$ in $X$ being the closure of $(\Delta, S(\cR))$. Let $(T_t)_{t\geq0}$ be the corresponding $C_0$-semigroup on $X$. Consider the  family of shift operators $(S_t)_{t\geq0}$,
\begin{align}\label{S(t)}
S_t\ffi(x):=\frac12\big(\ffi(x+\sqrt{t})+\ffi(x-\sqrt{t})\big),\qquad\forall\,\,\ffi\in X,\quad x\in\cR.
\end{align}
Then all $S_t$ are bounded linear operators on $X$, $\|S_t\|\leq1$ and for all $\ffi\in S(\cR)$ holds (via Taylor expansion):
\begin{align*}
S_t\ffi(x)-\ffi(x)&=\frac12\big(\ffi(x+\sqrt{t})-\ffi(x)\big)+\frac12\big(\ffi(x-\sqrt{t})-\ffi(x)\big)\\
&
=\frac12\left(\sqrt{t}\ffi'(x)+\frac12 t\ffi''(x)+o(t)  \right)+\frac12\left(-\sqrt{t}\ffi'(x)+\frac12t\ffi''(x)+o(t)  \right)\\
&
=t\ffi''(x)+o(t)=L\ffi(x)+o(t).
\end{align*}
Moreover, it holdts that $\lim_{t\to0}\|t^{-1}(S_t\ffi-\ffi)-L\ffi\|_X=0$ for all $\ffi\in S(\cRd)$. Hence the family 
$(S_t)_{t\geq0}$ is Chernoff equivalent to the heat semigroup~\eqref{heatSG} on $X$. Extending  $(S_t)_{t\geq0}$ to the $d-$dimensional case and applying the ``rotation'' techniques in $X=L^2(\cRd)$, one obtains Chernoff approximation for the Schr\"{o}dinger group $(e^{it\Delta})_{t\geq0}$ (\cite{MR3871530}). Further, one can apply the techniques of Sections~\ref{subsec_additive}---~\ref{subsec_subordinate}, to construct Chernoff approximations for Schr\"{o}dinger groups generated by more complicated  differential and pseudo-differential operators.

Let us now combine this techniques with averaging. Averaging is an extension of the classical Daletsky-Lie-Trotter  formula (see Sec.~\ref{subsec_additive}) for the case when the generator $(L,\Dom(L))$ in a Banach space $X$ is not just a finite sum of linear operators $L_k$, but an integral: 
$$
L:=\int_{\mathcal{E}} L_\eps d\mu(\eps),
$$ 
where $\mathcal{E}$ is a set and $\mu$ is a suitable probability measure on (a  $\sigma-$algebra of subsets of) $\mathcal{E}$, and $L_\eps$ are linear operators in $X$ for all $\eps\in\mathcal{E}$. It turns out that (under some additional assumptions) the family $(F(t))_{t\geq0}$,
$$
F(t)\ffi:=\int_\mathcal{E} e^{tL_\eps}\ffi\, d\mu(\eps),\qquad\qquad\ffi\in X,
$$
is Chernoff equivalent to the semigroup $(e^{tL})_{t\geq0}$ on $X$. Consider now $X=C_\infty(\cRd)$ or $X=L^p(\cRd)$, $p\in[1,\infty)$. Let us now generalize the family $(S_t)_{t\geq0}$ of~\eqref{S(t)} to the following family $(U_\mu(t))_{t\geq0}$:
consider the family $(S_\eps(t))_{t\geq0}$, $S_\eps(t)\ffi(x):=\ffi(x+\sqrt{t}\eps)$ for all $\ffi\in X$ and for a fixed $\eps\in\cRd$; define the family $(U_\mu(t))_{t\geq0}$ by
$$
U_\mu(t)\ffi(x):=\int_\cRd S_\eps(t)\ffi(x)\,d\mu(\eps)\equiv \int_\cRd \ffi(x+\eps\sqrt{t})\,d\mu(\eps).
$$
Assume that $\mu$ is a symmetric measure with finite (mixed) moments up to the third order and  positive second moments  $a_j:=\int_\cRd \eps_j^2\mu(d\eps)>0$, $j=1,\ldots,d$. Then one can show that the family $(U_\mu(t))_{t\geq0}$ is Chernoff equivalent to the heat semigroup $(e^{t\Delta_A})_{t\geq0}$, where $\Delta_A:=\frac12\sum_{j=1}^da_j\frac{\pd^2}{\pd x_j^2}$.  Substituting $(S_\eps(t))_{t\geq0}$ by the family $(S^\sigma_\eps(t))_{t\geq0}$, 
$S^\sigma_\eps(t)\ffi(x):=\ffi(x+\eps t^\sigma)$, for some suitable  $\sigma>0$, and choosing proper measures $\mu$, one can construct analogous Chernoff approximations for semigroups generated by fractional Laplacians and relativistic Hamiltonians. This approach can be further generalized by considering pseudomeasures $\mu$, what leads to Chernoff approximations for Schr\"{o}dinger groups.

\section*{Acknowledgements}
I would like to thank Christian Bender whose kind support made this work possible. I would like to thank Markus Kunze for  usefull discussions and communication of references~\cite{MR2541275,MR2812574}.

\bibliographystyle{abbrv}

\bibliography{Mabib_02-04-2019}

\end{document}